\documentclass[11pt]{article}
\usepackage{latexsym}
\usepackage{amsfonts,amssymb,amsmath,mathtools}
\usepackage{bbold}

\newtheorem{thm}{Theorem}[section]
\newtheorem{cor}[thm]{Corollary}
\newtheorem{lem}[thm]{Lemma}
\newtheorem{pro}[thm]{Proposition}
\newtheorem{defn}[thm]{Definition}

\newcommand{\ov }{\overline }
\newcommand{\minus}{\smallsetminus}
\newcommand{\x}{\hspace{-0.025in}\times\hspace{-0.025in}}
\newcommand{\e}{\varepsilon}
\title{A monoid version of the Brin-Higman-Thompson groups}

\author{J.C.\ Birget}

\date{\scriptsize 
27 vi 2020}

\begin{document}
\maketitle

\begin{abstract}
We generalize the Brin-Higman-Thompson groups $n G_{k,1}$ to monoids
$n M_{k,1}$, for $n \ge 1$ and $k \ge 2$, by replacing bijections by 
partial functions.
The monoid $n M_{k,1}$ has $n G_{k,1}$ as its group of units, and is 
congruence-simple. 
Moreover, $n M_{k,1}$ is finitely generated, and for $n \ge 2$ its word
problem is {\sf coNP}-complete.
We also present new results about higher-dimensional joinless codes.
\end{abstract}

%%%%%%%%%%%%%%%%%%%%%%%%%%%%%%%%%%%%%%%%%%%%%%%%%%%%%%%%%%%%%%%
% Sect. 1
%%%%%%%%%%%%%%%%%%%%%%%%%%%%%%%%%%%%%%%%%%%%%%%%%%%%%%%%%%%%%%%
\section{Introduction}

We extend the Brin-Higman-Thompson group $n G_{k,1}$ to a monoid, 
$n M_{k,1}$.  For this we follow in outline the definition of the monoid 
$M_{k,1}$ given in \cite{BiThompsMonV3, Equiv} (which extends the 
Thompson-Higman group $G_{k,1}$), and the string-based definition of the 
Brin-Thompson group $n G_{2,1}$ \cite{BinG} and the Brin-Higman-Thompson 
group $n G_{k,1}$ \cite{BinGk1}.

The groups $n G_{2,1}$ were introduced by Brin \cite{Brin1} for $n \ge 2$. 
For $n=1$, $\,1 G_{2,1} = G_{2,1}$ ($= V$) is the well-known Thompson group
\cite{Th0, CFP}. 
For $n \ge 2$, $n G_{2,1}$ ($= nV$) is an $n$-dimensional generalization of 
$V$.
Many interesting properties have been proved about $nG_{2,1}$: It is an 
infinite, finitely presented, simple group, and the groups $n G_{2,1}$ are 
non-isomorphic for different $n\,$  
\cite{Brin1, Brin2, Brin3, HennigMattucci, BleakLanoue}. 
In \cite{BinG}, a string-based description of $n G_{2,1}$ was introduced
(hinted at in \cite[Subsection 4.3]{Brin1}, and partly developed in 
\cite{BleakLanoue}). 
The Brin-Higman-Thompson group $n G_{k,1}$, for $n \ge 1$ and $k \ge 2$, is 
the obvious common generalization of the Brin-Thompson group $n G_{2,1}$ and 
the Higman-Thompson group $G_{k,1}$ \cite{Hig74}. In the literature 
$n G_{k,1}$ has been studied in a number of papers (one of the earliest is 
\cite{DicksMartinez}).
As $n G_{k,1}$ can be defined as a group of permutations of 
$n A_k^{\,\omega}$, where $A_k = \{a_0, a_1,\ldots, a_{k-1}\}\,$, it is 
natural to generalize this to a monoid $n M_{k,1}$ of (partial) functions on
$n A_k^{\,\omega}$.
When $n = 1$, $\,1 M_{k,1}$ is the monoid $M_{k,1}$ introduced in 
\cite{BiThompsMon, BiThompsMonV3}.

In the construction of $n M_{k,1}$ and the proofs of some of its 
properties, new results about {\em joinless codes} in $nA_k^{\,*}$ are 
proved.  In particular, we introduce the concepts of completion and 
complement of a non-maximal finite joinless code, and we show their 
existence. The definition of a right-ideal morphism leads to complications
that are not seen in the group case.

\medskip

\noindent The monoid $n M_{k,1}$ has interesting properties:

\smallskip

\noindent $\bullet$ \ $n M_{k,1}$ has $n G_{k,1}$ as group of units, and
is regular, $\cal J$-0 simple, and congruence-simple.

\smallskip

\noindent $\bullet$ \ $n M_{k,1}$ is finitely generated.

\smallskip

\noindent $\bullet$ \ For $n \ge 2$, the word problem of $n M_{k,1}$ over a 
finite generating set is {\sf coNP}-complete. 

\smallskip

\noindent $\bullet$ \ An important motivation for the study of $2 M_{2,1}$ 
is that it provides a way to describe acyclic digital circuits by words over 
a finite generating set, with at most polynomial length increase, and such
that the word problem of the monoid is polynomial-time equivalent to the 
equivalence problem for acyclic circuits. 
So $2 M_{2,1}$ can be viewed as an algebraic model of computation by 
acyclic circuits.

%\bigskip

%\bigskip

\newpage %%%%%%%%%%%%%%%%%%%%%%%%%%%%%%%%%%%%%%%%%%%%%%%%%%%%%%%%%%%%%%

\noindent {\bf Terminology and notation}

\smallskip

\noindent
-- ``Function'' means partial function $f$: $X \to Y$, with domain 
${\rm Dom}(f) \subseteq X$, and image ${\rm Im}(f) \subseteq Y$. 
Usually, the sets $X$ and $Y$ will be free monoids $A^*$, or Cantor spaces 
$A^{\omega}$, or their direct powers $nA^*$ or $nA^{\omega}$.
If ${\rm Dom}(f) = X$, then $f$ is called a {\it total function} on $X$.

\smallskip

\noindent
-- $f|_S$, \ the restriction of a function $f$ to a set $S$

\smallskip

\noindent
-- $A^*$, \ the free monoid freely generated by $A$, a.k.a.\ the set of all
   strings over the alphabet $A$

\smallskip

\noindent
-- $\e$, \ the empty string

\smallskip

\noindent
-- $|x|$, \ the length of the string $x \in A^*$

\smallskip

\noindent
-- $x \le_{\rm pref} y$, \ $x$ is a prefix of $y\,$
    (where $x \in A^*$ and $y \in A^* \cup A^{\omega}$)

\medskip

\noindent
-- $nA^*$  $=$ {\large \sf X}$_{_{i=1}}^{^n} A^*$, 
 \ \ $nA^{\omega}$  $=$  {\large \sf X}$_{_{i=1}}^{^n} A^{\omega}$, 
 \ \ the $n$-fold cartesian products

\medskip

\noindent
-- $x_i\,$ ($\,\in A^*$), \ the $i$th coordinate of $x \in nA^*\,$
  (for $1 \le i \le n$)

\smallskip

\noindent
-- ${\rm maxlen}(S) \,=\, \max\{|x_i|: x \in S, \ i \in \{1,\ldots,n\}\}$, 
 \ for $S \subseteq n A^*$.

\bigskip

%%%%%%%%%%%%%%%%%%%%%%%%%%%%%%%%%%%%%%%%%%%%%%%%%%%%%%%%%%%%%%%%%%%%%%
\section{Right ideals of {\boldmath $nA^*$}}

Let $A$ be an alphabet (in this paper, an alphabet is any finite non-empty 
set); typically we use the alphabets $A_1 = \{a_0\} \subset A_2 = \{a_0,a_1\}$
$\subseteq$ $A_k = \{a_0,a_1,\ldots,a_{k-1}\}$ for any $k \ge 2$; so 
$|A_k| = k$. For $x \in A^*$ the length of $x$ is denoted by $|x|$.
For $x_1, x_2 \in A^*$ the {\it concatenation} is denoted by $x_1 x_2$ or
$x_1 \cdot x_2$; it has length $|x_1| + |x_2|$.
For two subsets $S_1, S_2 \subseteq A^*$, we define the concatenation by
$\, S_1 \cdot S_2 = \{x_1 \cdot x_2 :$ $x_1 \in S_1$ and $x_2 \in S_2\}$.
For $x, p \in A^*$ we say that $p$ is a {\it prefix} of $x$ iff
$x = pu$ for some $u \in A^*$; this is denoted by $p \le_{\rm pref} x$.
Two strings $x, y \in A^*$ are called {\it prefix-comparable} (denoted
by $x \,\|_{\rm pref}\, y \,$) iff
$\, x \le_{\rm pref} y \,$ or $ \, y \le_{\rm pref} x$.
A {\it prefix code} (a.k.a.\ a prefix-free set) is any subset $P \subset A^*$
such that for all $p_1, p_2 \in P$: $\, p_1 \,\|_{\rm pref}\, p_2 \,$
implies $p_1 = p_2$; so a prefix code is the same thing as a 
$<_{\rm pref}$-antichain.

For $n \ge 1$, $n A^*$ denotes the $n$-fold
direct product {\large \sf X}$_{_{i=1}}^{^n} A^*$, with coordinate-wise
concatenation as multiplication.
The identity element of $n A^*$ is $(\e)^n$, i.e., the $n$-tuple 
consisting of $n$ copies of the empty string $\e$.
We also consider the $n$-fold cartesian product $n A^{\omega}$  $=$
{\large \sf X}$_{_{i=1}}^{^n} A^{\omega}\,$; this is the $n$-dimensional 
Cantor space. 
For a string $x \in A^{\ell}$ (with $\ell \in {\mathbb N}$), we call $\ell$ 
the {\em length} of $x$; we denote the length of $x$ by $|x|$.
For a non-empty finite set $S \subseteq n A^*$ we consider the {\em maximum 
length} in $S$, denoted by ${\rm maxlen}(S)$, and defined by:
$\, {\rm maxlen}(S)$  $\,=\,$
$\max\{|x_i| : x \in S, \ i \in \{1,\ldots,n\} \}$, where $x_i \in A^*$ is 
the $i$th coordinate of $x \in n A^*$.
In the case of a singleton $\{u\}$ we also write ${\rm maxlen}(u)$ for
${\rm maxlen}(\{u\})$.

The {\em initial factor order} is defined for $u, v \in nA^*$ by 
$u \le_{\rm init} v$ iff there exists $x \in nA^*$ such that $u x = v$.
Hence, $u \le_{\rm init} v$ in $nA^*$ iff
$u_i \le_{\rm pref} v_i$ in $A^*$ for all $i= 1, \, \ldots, n$.
The partial order $\le_{\rm init}$ on $nA^*$ is a generalization of the 
prefix order $\le_{\rm pref}$ on $A^*$; we chose a different name than 
``prefix order'' because for $n \ge 2$, $\le_{\rm init}$ has very different 
properties than $\le_{\rm pref}$ (as we shall see when we consider initial 
factor codes and joinless codes).
Note that $\le_{\rm init}$ is the same as $\ge_{\cal R}$, i.e., it is 
the reverse of the $\cal R$-order in the monoid $nA^*$.
An {\em initial factor code} is, by definition, a $<_{\rm init}$-antichain
(i.e., a set of elements of $nA^*$ no two of which are 
$<_{\rm init}$-related).

\begin{lem} \label{LEMprefEquOmega}  {\rm (folklore).} 
 \ For all strings $p, q \in A^*$ with $p \ne q$, there exists 
$u \in A^{\omega}$ such that $\,pu \ne qu$. 
\end{lem}
{\sc Proof.} In case $p$ and $q$ are not prefix-comparable, $p$ and $q$ 
have a longest common prefix $z \in A^*$ such that $|z| < |p|$ and 
$|z| < |q|$. So $p = zax$ and $q = zby$ for some $a, b \in A$ with 
$a \ne b$, and $x, y \in A^*$. It follows that for all $u \in A^{\omega}$, 
$\,p u = zax u \ne zbyu = q u$.

The remaining case is $p <_{\rm pref} q\,$ (the case $q <_{\rm pref} p$ is 
similar).  Then $p = q z$ for some $z \in A^*$ with $z \ne \e$. 
Let $y \in A^*$ be such that $y \ne z$ and $|y| = |z|$.
Then for $u = y^{\omega}$ we have $p u \ne q u$.
 \ \ \ $\Box$

\bigskip

\noindent The {\em join} of $u,v \in nA^*$ (with respect to the partial 
order $\le_{\rm init}$) is by definition  

\smallskip

 \ \ \  \ \ \  $u \vee v$   $\,=\,$
$\min_{\le_{\rm init}}\{z \in nA^* :$ $ \, u \le_{\rm init} z$ and 
$v \le_{\rm init} z\}$.

\smallskip

\noindent Of course, $u \vee v$ might not exist.
The following was proved in \cite[Lemma 2.5]{BinG}:

\begin{lem} \label{LEMjoinCharacteriz} {\bf (existence of joins).}
 \ For any $u,v \in nA^*$: 

The join $\,u \vee v \,$ exists \ iff 
 \ $\, u_i \,\|_{\rm pref}\, v_i \ $ for all $i \in \{1,$  $\ldots,$ 
$n\}.$ 

If $u \vee v$ exists then
$\, (u \vee v)_i \,=\, u_i\,$ if $v_i \le_{\rm pref} u_i$; 
 \ and $\, (u \vee v)_i \,=\, v_i \,$ if $u_i \le_{\rm pref} v_i$. So the
join is unique when it exists.
 \ \ \ $\Box$
\end{lem}
Lemma \ref{LEMjoinCharacteriz} immediately implies:

\begin{cor} \label{JoinVSinit}
 \ For all $u,v \in nA^*$: \ If $u \vee v$ exists, and if
$\,{\rm maxlen}(u) \le |v_i|\,$ for all $i \in$ $\{1,$  $\ldots,$ $n\}$, then 
$u \vee v = v$, hence $\,u \le_{\rm init} v$.
 \ \ \  \ \ \  $\Box$
\end{cor}

\begin{lem} \label{LEMjoinOmega}
 \ For all $x, x' \in nA^*$:  $ \ x \vee x'\,$ exists \ iff \ there are 
$u, u' \in n A^{\omega}$ such that $\,x u = x' u'$.  
\end{lem}
{\sc Proof.} Note that here, $u$ and $u'$ are in $n A^{\omega}$, not in
$nA^*$.   

The left-to-right implication is immediate from the definition 
of the join. Conversely, suppose $w = x u = x' u' \in n A^{\omega}$.
Then for every $i \in \{1,\ldots,n\}$: $\,x_i, x'_i \in A^*$ are prefixes of
$w_i \in A^{\omega}$. Hence $x_i$ and $x'_i$ are prefix-comparable.
Therefore, by Lemma \ref{LEMjoinCharacteriz}, $\,x \vee x'$ exists.
 \ \ \ $\Box$

\begin{lem} \label{LEMyNzNOJ}
 \ For every $y,z \in nA^*$ with $y \neq z$, there exists $v \in nA^*$ 
such that $\,yv \vee z\,$ does not exist, or $\,y \vee zv\,$ does not exist.
\end{lem}
{\sc Proof.} If $y \vee z$ does not exist then the Lemma holds with 
$v = (\e)^n$.
If $y \vee z$ exists and $y \neq z$, then Lemma \ref{LEMjoinCharacteriz} 
implies that $y_i = z_i a s >_{\rm pref} z_i$, or 
$z_i = y_i a s >_{\rm pref} y_i$, for some $i \in \{1,\ldots,n\}$, 
$\,a \in A$, and $\,s \in A^*$. Consider the first case (the other one 
being similar). In that case, for any $b \in A \minus \{a\}$: $\,y_i b$ is 
not prefix-comparable with $z_i = y_i a s$. Hence (by Lemma
\ref{LEMjoinCharacteriz}): $\, y \cdot v \vee z\,$ does not exist, where 
$v = (\e)^{i-1} \times (b) \times (\e)^{n-i}$.
  \ \ \ $\Box$

\bigskip

By definition, an {\em initial factor code} is a set $S \subseteq nA^*$ such
that no two different elements of $S$ are $<_{\rm init}$-comparable.
A {\em maximal initial factor code} is an initial factor code $S$ that is
not a strict subset of any other initial factor code of $nA^*$.

A {\em joinless code} in $n A^*$ is a subseteq $S$ such that no two different 
elements have a join.  
A {\em maximal joinless code} in $nA^*$ is a joinless code 
$C \subseteq n A^*$ such that $C$ is not a strict subset of any other 
joinless code in $nA^*$.
In the geometric interpretation of $nA^*$, a maximal joinless code is a
{\em tiling} of the hypercube $[0,1]^n\,$ (see \cite{BinG} for examples and
more information on these definitions).

In this paper we only consider {\em finite} joinless codes and initial 
factor codes.

\begin{defn} \label{DEFtiling} {\bf (maximal joinless code within a set).}
 \ Let $u \in nA^*$.  A joinless code $C \subseteq nA^*$ is {\em maximal in}
$u \cdot nA^*$ \ iff \ $C \subseteq u \cdot nA^*$, and $C$ is not a strict 
subset of any other joinless code in $\,u \cdot nA^*$. 
Then we also say that $C$ is a {\em tiling of} $u$.
\end{defn}
Geometrically, a joinless code $C$ is maximal in $u \cdot nA^*\,$ iff
$\,C$ is a tiling of the hyperrectangle represented by $u$.
This is the case iff $\,C = u \cdot Q$ for some maximal joinless code 
$Q \subseteq n A^*$. 

\bigskip

A right ideal of $n A^*$ is defined as in any semigroup: $R \subseteq n A^*$
is a right ideal iff $R \cdot n A^* \subseteq R$. 
A generating set of a right ideal $R$ is a subset $C \subseteq R$ such that
$R = C \cdot n A^*$. 

For every right ideal $R$ there exists a unique maximal initial-factor code 
that generates $R\,$ \cite[Lemma 2.7]{BinG}.

\begin{defn} \label{DEFsetJoin} {\bf (join of sets).}
 \ For any sets $X, Y \subseteq n A^*$, the {\em join} is defined by

\smallskip
 
 \ \ \ $X \vee Y \ = \ \{x \vee y \,\in\, n A^* \,:\, x \in X, \ y \in Y\}$
\end{defn}
See \cite[Prop.\ 2.18]{BinG}.

For all $X \subseteq n A^*$: $X \subseteq X \vee X$. If $X$ is joinless 
then $X = X \vee X$; but in general, $X = X \vee X$ does not imply that 
$X$ is joinless (a counter-example is $X = \{a_0, a_0 a_0\}$ with $n=1$ and 
$A = \{a_0, a_1\}$).

\begin{lem} \label{joinAssoc} {\bf (associativity of join).}
 \ For all $u,v,w \in n A^*$, statements {\rm (1) - (4)} are equivalent:

\smallskip

\noindent {\small \rm (1)} \ \ $u \vee (v \vee w)$ \ exists;

\smallskip

\noindent {\small \rm (2)} \ \ $(u \vee v) \vee w$ \ exists;

\smallskip

\noindent {\small \rm (3)} \ \ for all $i \in \{1, \ldots,n\}:$ \ the three
strings $u_i$, $v_i$, and $w_i$ are prefix-comparable two-by-two;

\smallskip

\noindent {\small \rm (4)} \ \ there exists a common upper-bound for 
$u, v, w$ with respect to $\le_{\rm init}$.

\medskip

\noindent Moreover:

\smallskip

\noindent {\rm (a)} \ If {\small \rm (1)} holds then
$\, u \vee (v \vee w) = (u \vee v) \vee w\,$;
 \ the common value is denoted by $\,u \vee v \vee w$.

 \ And for all $i:$ \ \ $(u \vee v \vee w)_i$  $\,=\,$
$u_i \vee_{\rm pref} v_i \vee_{\rm pref} w_i$.

\smallskip

\noindent {\rm (b)} \ For all sets $X, Y, Z \subseteq n A^*${\rm :}
 \ \ $X \vee (Y \vee Z) = (X \vee Y) \vee Z$.
\end{lem}
{\sc Proof.} (1) - (4) and (a) are straightforward consequences of Lemma
\ref{LEMjoinCharacteriz}.

(b) Note that the join of sets always exists (it may be the empty set).
By definition of the join of sets, $(X \vee Y) \vee Z$  $=$
$\{(x \vee y) \vee z: x \in X,\, y \in Y,\, z \in Z\}$. By item (a), this is
equal to $\{x \vee (y \vee z): x \in X,\, y \in Y,\, z \in Z\}$  $=$
$X \vee (Y \vee Z)$.
 \ \ \ $\Box$

%\medskip
%
%By Lemma \ref{joinAssoc} we can consider multiple joins of the form
% $\,x^{(1)} \vee \,\ldots\, \vee x^{(k)}\,$ (for $k \ge 1$). 
% A join of $k$ two-by-two different elements of $n A^*$ is called a 
% {\em join of multiplicity} $k$.

\begin{lem} \label{joinDistrCup} {\bf (distributivity of 
{\boldmath $\vee$ over $\cup$}).}
 \ For all sets $X, Y, Z \subseteq n A^*:$

\smallskip

$X \vee (Y \cup Z) \,=\, (X \vee Y) \,\cup\, (X \vee Z)$, \ \ and
 \ \ $(X \cup Y) \vee Z \,=\, (X \vee Z) \,\cup\, (Y \vee Z)$.
\end{lem}
{\sc Proof.} By definition of the join of sets, $X \vee (Y \cup Z)$  $=$
$\{x \vee u : x \in X,\, u \in Y \cup Z\}$. By the definition of $\cup$
this is equal to $ \ \{x \vee u : x \in X,\, u \in Y\}$   $\,\cup\,$
$\{x \vee u : x \in X,\, u \in Z\}$  $=$
$(X \vee Y) \,\cup\, (X \vee Z)$.

Right distributivity is proved in a similar way.
 \ \ \ $\Box$

\medskip

\noindent The main importance of the set join is its connection with the
intersection of right ideals:

\begin{lem} \label{SetJoin} {\bf ({\boldmath $\vee$ of sets, and $\cap$ of
ideals}). }
 \ For any sets $X, Y \subseteq n A^*$:

\smallskip

\noindent {\small \rm (1)} \ \ $X \cdot n A^* \,\cap\, Y \cdot n A^*$  
$ \ = \ $   $(X \vee Y) \cdot n A^*$, \ and

\smallskip

 \ \ $X \cdot n A^{\omega} \,\cap\, Y \cdot n A^{\omega}$  $ \ = \ $
$(X \vee Y) \cdot n A^{\omega}$.

\smallskip

\noindent {\small \rm (2)} \ \ If $X$ and $Y$ are finite then
$\, X \cdot n A^* \,\cap\, Y \cdot n A^*\,$ is finitely generated (as a
right ideal).

\smallskip

\noindent {\small \rm (3)} \ \ If $X$ and $Y$ are joinless codes then
$\, X \cdot n A^* \,\cap\, Y \cdot n A^*\,$ is joinless generated (by the
joinless code $X \vee Y$).
\end{lem}
{\sc Proof.} For (1) in $n A^*$, see the proof of
\cite[Prop.\ 2.18(3)]{BinG}.  (Note that \cite[Prop.\ 2.18(3)]{BinG} 
assumes that $X$ and $Y$ are joinless; but the proof actually does not 
depend on that.)

Let us prove (1) for $n A^{\omega}$.
For all $w \in X \cdot n A^{\omega} \,\cap\, Y \cdot n A^{\omega}$:  
$w$ has an initial factor $x \in X$ and an initial factor $y \in Y$. Then
$x \vee y$ is also an initial factor of $w$ (by Lemma \ref{LEMjoinCharacteriz}). So, $w \in (X \vee Y) \cdot n A^{\omega}$.
Conversely, if $w \in (X \vee Y) \cdot n A^{\omega}$, then $w$ has an initial
factor $x \vee y \in X \vee Y$. Since $x$ and $y$ are initial factors of
$x \vee y$, it follows that $w$ has an initial factor in $X$ and an initial
factor in $Y$; hence $w \in X \cdot n A^* \,\cap\, Y \cdot n A^*$.

(2) is an immediate consequence of (1) and the definition of $X \vee Y$.

(3) is proved in \cite[Prop.\ 2.18(1)]{BinG}.
 \ \ \   \ \ \ $\Box$

\begin{lem} \label{LEMdisjointOmeaga}
 \ For all finite joinless codes $P, Q \subseteq n A^*$ the following 
are equivalent:

\smallskip

\noindent {\small \rm (1)} \ \ \ $P \vee Q = \varnothing \,$;

\smallskip

\noindent {\small \rm (2)} 
 \ \ \ $P \cdot n A^* \ \cap \ Q \cdot n A^* \,=\, \varnothing\,$;
 
\smallskip

\noindent {\small \rm (3)} 
 \ \ \ $P \cdot n A^{\omega} \ \cap \ Q \cdot n A^{\omega}$  $\,=\,$   
$\varnothing\,$.
\end{lem}
{\sc Proof.} This follows immediately from Lemma \ref{SetJoin}.
 \ \ \ $\Box$

\bigskip

It is easy to prove that for all $P, Q \subseteq n A^*$: 
$ \, P \cdot n A^* = Q \cdot n A^*\,$ implies
$P \cdot n A^{\omega}$ $=$  $Q \cdot n A^{\omega}$.
However, {\em the converse does not hold}. 
E.g., for $n=1$, $A = \{a_0, a_1\}$, $\,P = \{a_0, a_1\}$, and 
$\,Q = \{\e\}$, we have 
$P A^{\omega} = Q A^{\omega}$ $=$ $A^{\omega}$, but $P A^* \ne Q A^*$. 
The following Definition and Lemma characterize when a pair of finite sets 
$P, Q \subseteq n A^*$ satisfies $P \cdot n A^{\omega}$ $=$  
$Q \cdot n A^{\omega}$.

\begin{defn} \label{DEFequivFin} {\bf (equivalences  {\boldmath
$\,\equiv_{\rm fin}$ and $\,\equiv_{\rm bd}$}).} 
 \ Let $P, Q \subseteq n A^*$ be any sets.

\smallskip

\noindent {\smallskip \rm (1)} \ The relation $\equiv_{\rm fin}$ between 
sets is defined by

\smallskip

 \ \ \   \ \ \ $P \equiv_{\rm fin} Q$ \ \  iff
 \ \ $P \cdot n A^* \,\vartriangle\, Q \cdot n A^*$ \ is finite,

\smallskip

 \ where $\vartriangle$ denotes symmetric difference. 

\smallskip

\noindent {\smallskip \rm (2)} \ The relation $\equiv_{\rm bd}$ between 
sets is defined by

\smallskip

 \ \ \   \ \ \ $P \equiv_{\rm bd} Q$ \ \  iff
 \ \ $P \cdot n A^{\omega} = Q \cdot n A^{\omega}$.
\end{defn}
In this paper we will use these equivalences only between finite sets.

The relations $\equiv_{\rm bd}$ and $\equiv_{\rm end}$ (defined below) were
introduced in \cite{Equiv}; the subscript {\sf bd} stands for {\em bounded 
end-equivalence}, and the subscript {\sf end} stands for 
{\em end-equivalence}.

\begin{lem} \label{LEMequivFinOmega} {\bf ({\boldmath $\equiv_{\rm fin}$ 
and $n A^{\omega}$}).} 
 \ For any {\em finite} sets $P, Q \subseteq n A^*$ the following are
equivalent:

\smallskip

\noindent {\small \rm (1)} \ \ $P \equiv_{\rm fin} Q$;

\smallskip

\noindent {\small \rm (2)} \ \ $P \equiv_{\rm bd} Q$;

\smallskip

\noindent {\small \rm (3)}
 \ \ both 
 \ $(\forall p \in P)(\forall z \in nA^*)(\exists q \in Q)[\,pz \vee q\,$
${\rm exists}\,]$ \ and 

 \ \ $\,(\forall q \in Q)(\forall z \in nA^*)(\exists p \in P)[\,p \vee qz\,$
${\rm exists}\,]$ \ \ hold.
\end{lem}
{\sc Proof.} $[(1) \Rightarrow (2)] \,$ If
$P \cdot n A^* \vartriangle Q \cdot n A^*\,$ is finite, let $\,n_0 = $
${\rm maxlen}(P \cdot n A^* \vartriangle Q \cdot n A^*)$.
Then for $z \in P \cdot n A^*$ with ${\rm maxlen}(z) > n_0$:
 \ $z \in Q \cdot n A^*$.
Indeed, if we had $z \not\in Q \cdot n A^*$, then
$z \in P \cdot n A^* \vartriangle Q \cdot n A^*$; but since
${\rm maxlen}(z) > n_0$, this would contradict the definition of $n_0$.
Similarly, if $z \in Q \cdot n A^*$ and ${\rm maxlen}(z) > n_0$, then
$z \in P \cdot n A^*$.
Thus for all $z \in nA^*$ with ${\rm maxlen}(z) > n_0$:
$z \in P \cdot n A^*$ iff $z \in Q \cdot n A^*$; in other words,
$P \cdot n A^* \,\cap\, n A^{\ge n_0}$  $=$
$Q \cdot n A^* \,\cap\, n A^{\ge n_0}$.
Item (2) then follows.

$[(2) \Rightarrow (3)] \,$ Since $P \cdot n A^{\omega}$ $\subseteq$
$Q \cdot n A^{\omega}$, every $p z \in P \cdot nA^*$ is an initial factor of
some $q w$ with $q \in Q$ and $w \in n A^{\omega}$. Therefore
every $p_i z_i$ is prefix-comparable with $q_i$, for all $i$. By 
\cite[Lemma 2.5]{BinG}, this implies that $p z \vee q$ exists.
% In particular, for $z = (\e)^n$, $p \vee q$ exists.

From $Q \cdot n A^{\omega}$ $\subseteq$ $P \cdot n A^{\omega}$  we derive
the second clause.

$[(3) \Rightarrow (1)] \,$ Let $\, n_0 = {\rm maxlen}(P \vee Q)$, which
exists since $P$ and $Q$ are finite.
If $p z \in P \cdot nA^*$ is such that ${\rm maxlen}(pz) > n_0$, then
$pz \in (P \vee Q) \cdot nA^*$ ($= P \cdot nA^* \,\cap\, Q \cdot nA^*$, by
Lemma \ref{SetJoin}).
Indeed, $pz \in P \cdot nA^*$, and by assumption, $pz \vee q$ exists for some
$q \in Q$; so $p_i z_i \,\|_{\rm pref}\, q_i$.
Since ${\rm maxlen}(pz) > n_0$, it follows that $p_i z_i >_{\rm pref} q_i$,
for all $i$.  Hence $pz \ge_{\rm init} q$. So
$pz \in P \cdot nA^* \,\cap\, Q \cdot nA^*$.
Therefore, if ${\rm maxlen}(pz) > n_0$ then
$pz \not\in P \cdot nA^* \,\vartriangle\, Q \cdot nA^*$. So, by
contraposition, if $pz \in P \cdot nA^* \,\vartriangle\, Q \cdot nA^*$,
then ${\rm maxlen}(pz) \le n_0$.

Similarly one proves for $qz \in Q \cdot nA^*$: if $qz \in $
$P \cdot nA^* \,\vartriangle\, Q \cdot nA^*$,
then ${\rm maxlen}(qz) \le n_0$.

This implies that $P \cdot nA^* \,\vartriangle\, Q \cdot nA^*$ is finite.
 \ \ \ $\Box$

\begin{lem} \label{LEMequivFTrans}
 \ The relation $\equiv_{\rm fin}$ on finite sets is transitive.
\end{lem}
{\sc Proof.}  This follows immediately from Lemma
\ref{LEMequivFinOmega}(1)(2).
 \ \ \ $\Box$

\begin{defn} \label{DEFendequiv} {\bf (end-equivalence).}
 \ We define the following relation between sets
$P, Q \subseteq nA^*$:

\smallskip

 \ \ \ $P \equiv_{\rm end} Q$ \ \ \ iff \ \ \ $(\forall x \in nA^*)$
$[\,P \cdot nA^* \,\cap\, x \cdot nA^* = \varnothing$ \ \ $\Leftrightarrow$
 \ \ $Q \cdot nA^* \,\cap\, x \cdot nA^* = \varnothing \,]$.
\end{defn}
In other words, $P \equiv_{\rm end} Q$ \ iff \ $P \cdot nA^*$ and 
$Q \cdot nA^*$ intersect the same right ideals of $nA^*$.

\bigskip

\noindent  {\bf Remarks on the topology of {\boldmath $n A^{\omega}$}:} 
The topology of the
$n$-dimensional Cantor space $n A^{\omega}$ is determined by
$\,\{L \cdot n A^{\omega} : L \subseteq A^*\}$, as the set of open sets.
%% We will denote this topological space by $n {\frak C}$.
The closure of a set $S \subset n A^{\omega}$ is denoted by ${\sf cl}(S)$,
and the interior by ${\sf in}(S)$.  The following is not hard to prove. 
For all $P \cdot n A^{\omega}$ and $Q \cdot n A^{\omega}$, where 
$P, Q \subset n A^*$:

\smallskip

 \ \ \ ${\sf cl}(P \cdot n A^{\omega}) = {\sf cl}(Q \cdot n A^{\omega})$ 
 \ \ \ iff \ \ \ $P \equiv_{\rm end} Q$.

\smallskip

\noindent See \cite{Equiv}, where this was considered for $A^*$.
It is easy to prove that 

\smallskip

 \ \ \ $(\forall x \in n A^*)$
$[\,P \cdot nA^{\omega} \cap x \cdot nA^{\omega} = \varnothing$
$ \ \Leftrightarrow \ $  
$Q \cdot nA^{\omega} \cap x \cdot nA^{\omega} = \varnothing\,]$

\smallskip

 \ \ \ iff \ \ \ $(\forall x \in n A^*)$
$[\,P \cdot nA^* \cap x \cdot nA^*  = \varnothing$
$ \ \Leftrightarrow \ $  
$Q \cdot nA^* \cap x \cdot nA^* = \varnothing\,]$.

\smallskip

\noindent When $P$ and $Q$ are {\em finite} then by Lemmas 
\ref{LEMequivFinOmega} and \ref{LEMequivFinEnds}: 

\smallskip

 \ \ \ $P \equiv_{\rm fin} Q$ \ iff
 \ $P \cdot nA^{\omega} = Q \cdot nA^{\omega}$ 
 \ iff \ $P \equiv_{\rm end} Q$.

\smallskip

\noindent In \cite{Equiv} it was proved in the case of $n=1$, that for 
infinite sets the relations  $\equiv_{\rm fin}$, $\equiv_{\rm bd}$, and 
$\equiv_{\rm end}$ are different. 
 \ \ \ [End, Remark.]

\begin{lem} \label{LEMequivFinEnds}
 \ For every {\em finite} sets $P, Q \subseteq nA^*$:
 \ \ \ $P \equiv_{\rm fin} Q$ \ \ iff \ \ $P \equiv_{\rm end} Q$ .
\end{lem}
{\sc Proof.} We use the fact that $P \equiv_{\rm fin} Q\,$ iff
$\,P \cdot n A^{\omega}$  $=$  $Q \cdot n A^{\omega}$ \ (Lemma
\ref{LEMequivFinOmega}).

\smallskip

\noindent $[\Rightarrow]$ Suppose $P \cdot n A^{\omega}$  $=$
$Q \cdot n A^{\omega}$. If $\,P \cdot nA^* \,\cap\, x \cdot nA^*$  $\ne$
$\varnothing$, then $P \cdot nA^{\omega} \,\cap\, x \cdot nA^{\omega}$
$\ne$  $\varnothing$, so $\,xw \in P \cdot nA^{\omega}$ for some $w \in$
$nA^{\omega}$. Hence $xw \in Q \cdot nA^{\omega}$
($\, = P \cdot nA^{\omega}$).
Therefore, $x$ is an initial factor of some string $q u$ for some $q \in Q$
and $u \in nA^*$.
Hence $Q \cdot nA^* \,\cap\, x \cdot nA^* = \varnothing$.

\smallskip

\noindent $[\Leftarrow]$ Suppose $\,P \cdot nA^* \cap x \cdot nA^*$  $=$
$\varnothing$ \ $\Leftrightarrow$
 \ $Q \cdot nA^* \cap x \cdot nA^* = \varnothing$.
Note that this means that in the $n$-dimensional Cantor space topology,
${\sf cl}(P \cdot nA^{\omega})$  $=$  ${\sf cl}(Q \cdot nA^{\omega})$,
where ${\sf cl}(.)$ denotes closure.
When $P$ is finite then $P \cdot nA^{\omega}$ is closed in the
$n$-dimensional Cantor space topology, i.e., ${\sf cl}(P \cdot nA^{\omega})$
$=$  $P \cdot nA^{\omega}$.  Similarly, $Q \cdot nA^{\omega}$ is closed.
Hence, $P \cdot n A^{\omega}$  $=$  $Q \cdot n A^{\omega}$.
 \ \ \ $\Box$

\begin{lem} \label{LEMmaxOmega} 
 \ For every {\em finite} joinless code $P \subseteq n A^*$ the following
are equivalent:

\smallskip

\noindent {\small \rm (1)} \ \ \ $P\,$ is maximal (as a joinless code);  
 
\smallskip

\noindent {\small \rm (2)} \ \ \ $P \cdot n A^{\omega} = n A^{\omega}\,$;  
 
\smallskip

\noindent {\small \rm (3)} \ \ \ $P \equiv_{\rm fin} \{\e\}^n\,$.
\end{lem}
{\sc Proof.} The equivalence of the last two statements follows immediately 
from Lemma \ref{LEMequivFinOmega}(1)(2).  
For the equivalence $(1) \Leftrightarrow (2)$, we will prove the 
contrapositive: $ \ P \cdot n A^{\omega} \ne n A^{\omega}$ \ iff
 \ $P$ is a not maximal. 

\smallskip

\noindent
$[\Leftarrow]$ By definition of maximality of a joinless code, 
$P$ is not maximal iff there exist $q \in n A^*$ such that $P \cup \{q\}$ is 
joinless. This holds iff 
$\,q \cdot nA^{\omega} \,\cap\, P \cdot nA^{\omega}$  $=$  $\varnothing$, 
which implies $\, P \cdot nA^{\omega} \ne n A^{\omega}$.

\smallskip

\noindent
$[\Rightarrow]$ If $P \cdot nA^{\omega} \ne n A^{\omega}$ then there exists 
$w \in n A^{\omega} \minus P \cdot nA^{\omega}$. 
Let $q \in nA^*$ be any initial factor of $w$ such that 
$\, |q_i| \ge {\rm maxlen}(P)\,$ for all $i \in \{1,\ldots,n\}$. 
Since $P$ is finite, ${\rm maxlen}(P)$ exists, and since every coordinate of
$w$ is infinitely long, $q$ exists.  

If, by contradiction, $q$ has a join with some $p \in P$ then (by Coroll.\ 
\ref{JoinVSinit}), $q \vee p = q$, so $p \le_{\rm init} q$. This implies
that $w \in p \cdot n A^{\omega}$, contrary to the assumption that 
$w \in n A^{\omega} \minus P \cdot nA^{\omega}$. Hence, since $q$ has no join
with any element of $P$, $P$ is not maximal.
 \ \ \ $\Box$

\begin{lem} \label{LEMequivFinIntersect}
 \ For any finite sets $P, Q \subseteq n A^*$: \ If $P \equiv_{\rm fin} Q\,$ 
then $ \ P \,\equiv_{\rm fin}\, Q \,\equiv_{\rm fin}\, P \vee Q$.
\end{lem}
{\sc Proof.} In this proof let us abbreviate $P \cdot n A^*$ by $P_A$, and
$Q \cdot n A^*$ by $Q_A$.  Let us show that 
$P_A \vartriangle (P \vee Q) \cdot n A^*$ is finite.
By Lemma \ref{DEFsetJoin}(1), $P_A \vartriangle (P \vee Q) \cdot n A^*$  
$=$  $P_A \vartriangle (P_A \cap Q_A)$. 

\smallskip

And $\, P_A \vartriangle (P_A \cap Q_A)$  $=$ 
$(P_A \,\cap\, \ov{P_A \cap Q_A}\,) \,\cup\, (\, \ov{P}_A \cap P_A \cap Q_A)$
$=$ $P_A \cap (\, \ov{P}_A \cup \ov{Q}_A \,)$
$=$ $P_A \cap \ov{Q}_A$ $=$ $P_A \minus Q_A$ $\subseteq$ 
$P_A \vartriangle Q_A$, which is finite. 

In a similar way one proves that $Q \vartriangle (P_A \cap Q_A)$ is finite. 
 \ \ \ $\Box$

\begin{lem} \label{LEMequivFPPnAell}
 \ For any finite set $S \subseteq n A^*$ and any $\ell \ge 0:$ 
 \ \ $S \,\equiv_{\rm fin}\, S \vee n A^{\ell}$.
\end{lem}
{\sc Proof.} By Lemma \ref{SetJoin}, 
$\,S \ n A^* \vartriangle (S \vee n A^{\ell}) \cdot n A^*$  $\,=\,$
$S \ n A^* \,\vartriangle\, (S \ n A^* \,\cap\, n A^{\ge \ell})$.
Since $S \ n A^* \,\cap\, n A^{\ge \ell} \,\subseteq\, S \ n A^*$, the 
latter symmetric difference is equal to 
$\, S \ n A^* \minus (S \ n A^* \,\cap\, n A^{\ge \ell})$, which is
equal to $\,\{x \in S \ n A^* \,:\, |x_i| < \ell$ for $i = 1,\ldots,n\}$.
This set is finite.
 \ \ \ $\Box$

\begin{lem} \label{LEMequivEqual} 
 \ For any $\ell \ge 0:$ \ If $\,S, T \subseteq n A^{\ell}$ and 
$S \equiv_{\rm fin} T$, then $S = T$.
\end{lem}
{\sc Proof.} Since $S \subseteq n A^{\ell}$, 
$\, {\sf init}(S \, n A^{\omega}) \,\cap\, n A^{\ell} = S$, where 
${\sf init}(S \, n A^{\omega})$ $ \ (\subseteq  n A^*)$ is the set of initial
factors of the elements of $S \, n A^{\omega}$.  Similarly, 
$\, {\sf init}(T \, n A^{\omega}) \,\cap\, n A^{\ell} = T$.
Since $S \equiv_{\rm fin} T$, we have 
$S \, n A^{\omega} = T \, n A^{\omega}$, by Lemma \ref{LEMequivFinOmega}. 
Hence, $S = T$.
 \ \ \ $\Box$

\begin{lem} \label{LEMequivFnAell}
 \ For any finite sets $P, Q \subseteq n A^*$ and any 
$\, \ell \ge {\rm maxlen}(P \cup Q):$

\smallskip

 \ \ \  \ \ \ $P \equiv_{\rm fin} Q$ \ \ iff 
 \ \ $P \vee n A^{\ell} \,=\, Q \vee n A^{\ell}$.
\end{lem}
{\sc Proof.} $[\Leftarrow]$ By Lemma \ref{LEMequivFPPnAell},
$P \equiv_{\rm fin} P \vee n A^{\ell}$ $\,=\,$  
$Q \vee n A^{\ell} \equiv_{\rm fin} Q$.

\noindent $[\Rightarrow]$ Suppose $P \equiv_{\rm fin} Q$, hence by Lemma
\ref{LEMequivFPPnAell}, 
$P \vee n A^{\ell} \equiv_{\rm fin} Q \vee n A^{\ell}$.  In addition, 
$\ell \ge {\rm maxlen}(P \cup Q)$ implies that 
$P \vee n A^{\ell}, \ Q \vee n A^{\ell} \subseteq n A^{\ell}$.
Now, Lemma \ref{LEMequivEqual} implies that 
$P \vee n A^{\ell} \,=\, Q \vee n A^{\ell}$.
 \ \ \ $\Box$

\begin{lem} \label{OneStepRestrjoinless} {\bf (one-step restriction and
extension).} 
 \ Let $P \subseteq n A^*$ be a finite set. For any
$\, p = (p_1, \, \ldots, p_n) \in P$ and $i \in \{1, \, \ldots, n\}$, let

\medskip

\hspace{0.3in} $P_{p,i} \ = \ (P \minus \{p\})$  $ \ \cup \ $
$\{(p_1, \,\ldots, p_{i-1},\, p_i a,\, p_{i+1},\, \ldots, p_n)\, :\,$
$a \in A\}$

\smallskip

\hspace{0.6in} $ \ = \  (P \minus \{p\})$ $ \ \cup \ $ 
$p \cdot (\{\e\}^{i-1} \x A \x \{\e\}^{n-i})$.

\medskip

\noindent Then we have:

\smallskip

\noindent {\small \rm (1)} \ $P$ is joinless \ iff \ $P_{p,i}\,$ is joinless.

\smallskip

\noindent {\small \rm (2)} \ $P \,\equiv_{\rm fin}\, P_{p,i}$.

\smallskip

\noindent {\small \rm (3)} \ $P$ is a maximal joinless code \ iff 
 \ $P_{p,i}\,$ is a maximal joinless code.

\medskip

\noindent The set $P_{p,i}$ is called a {\em one-step restriction} of $P$
(``restriction'' because
$\, P_{p,i} \cdot (nA^*) \,\subseteq\, P \cdot (nA^*)$);
and $P$ is called a {\em one-step extension} of $P_{p,i}$.
\end{lem}
{\sc Proof.} (1) This is proved in \cite[Lemma 2.11(1)]{BinG}. \\  
(2) Either $P \cdot n A^* = P_{p,i} \cdot n A^*$; or
$P \cdot n A^* \minus  P_{p,i} \cdot n A^* = \{p\}$,
and $\, P_{p,i} \cdot n A^* \subseteq P \cdot n A^*$. Hence,
$P \cdot n A^* \vartriangle P_{p,i} \cdot n A^*$ is either empty or
$\{p\}$, so it is finite. \\   
(3) This is proved in \cite[Lemma 2.11(2)]{BinG}; by Lemma
\ref{LEMmaxOmega}(3) it also follows from item (2) of the present Lemma.
 \ \ \ $\Box$

\medskip

\noindent {\bf Remark.} Lemma \ref{OneStepRestrjoinless}(1) applies to 
joinless codes, but it does not hold in a similar way for initial factor 
codes.  For example, consider the initial factor code
$\, P = \{(\e,a_0),\, (a_0,\e)\}$ in $\,2\,\{a_0,a_1\}^*$.
Then for $p = (a_0, \e)$ and $i = 2$ we obtain
$\, P_{p,i}' = \{(\e,a_0),\, (0_, a_0),\, (a_0, a_1)\}$, which is not an
initial factor code.

\begin{lem} \label{LEMonestepCartProd}
 \ {\small \rm (1)} \ If a set $P \subset nA^*$ is obtained from 
$\{\e\}^n$ by a finite sequence of one-step restrictions, then 
$P$ is a finite maximal joinless code.

The converse is true for $n \le 2$. However, for every $n \ge 3$ there exists
a finite maximal joinless code in $n A^*$ that cannot be obtained from 
$\{\e\}^n$ by a finite sequence of one-step restrictions.

\smallskip

{\small \rm (2)} \ If $P \subset mA^*$ and $Q \subset nA^*$ are finite 
maximal joinless codes that can be obtained from $\{\e\}^m$, 
respectively $\{\e\}^n$, by one-step restrictions, then 
$P \times Q \subset (m+n)A^*$ is a finite maximal joinless code that can be 
obtained from $\{\e\}^{m+n}$ by one-step restrictions. 
\end{lem}
{\sc Proof.} {\small (1)} Every one-step restriction preserves joinlessness 
and maximality (by Lemma \ref{OneStepRestrjoinless}(1)(3)). The converse for
$n=1$ is folklore. For $n = 2$ it was first proved in 
\cite[Thm.\ 12.11]{LawsonVdovina}, in a different formulation; in the case 
of $A = \{a_0, a_1\}$ a different proof is given in 
\cite[Lemma 2.10(2)]{BinG}.  
For $n \ge 3$, a counter-example to the converse appears in 
\cite[Ex.\ 12.8]{LawsonVdovina}.

{\small (2)} This follows from the fact that a one-step reduction is 
applied to one coordinate, independently of the other coordinates. Let us 
denote existence of a sequence of one-step reductions from a set $X$ to a 
set $Y$ by $X \stackrel{*}{\to} Y$. Then 
 \ $\{\e\}^{m+n}$ 
$ \ \stackrel{*}{\to} \ $   $P \times \{\e\}^n$ 
$ \ \stackrel{*}{\to} \ $   $P \times Q$.
 \ \ \ $\Box$

\begin{lem} \label{LEMequivStoSnAL}
 \ For any finite set $S \subseteq n A^*$ and for any integer
$\,\ell \ge {\rm maxlen}(S):$   $ \ S \vee n A^{\ell}$ can be reached from 
$S$ by a finite sequence of one-step restrictions.
\end{lem}
{\sc Proof.} \ We have 

\smallskip

 \ \ \  \ \ \ $S \vee n A^{\ell}$  $ \ = \ $
$\big(\bigcup_{s \in S} \{s\} \big) \vee n A^{\ell}$  $ \ = \ $ 
$\bigcup_{s \in S} (\{s\} \vee n A^{\ell})$, 

\smallskip

\noindent the latter by distributivity (Lemma \ref{joinDistrCup}).  So it 
suffices now to show that $\{s\} \vee n A^{\ell}$ is reached from $\{s\}$ 
by one-step restrictions.  We prove this by induction on 
$\, N(s) \,=\, n \ell \,-\, \sum_{i=1}^n |s_i|$.  

\smallskip 

If $N(s) =0$ then $\ell = {\rm maxlen}(S)$ and all the coordinates of $s$ 
have length $\ell$, so $s \in n A^{\ell}$, hence 
$\, \{s\} \vee n A^{\ell} = \{s\}$.  Of course, $\{s\}$ is reachable from 
$\{s\}$ (in 0 steps). 

If $N(s) > 0$, then $|s_i| < \ell$ for some $i$. Applying a one-step 
restriction to $s$ at such a coordinate $i$ yields
$\, \{(s_1, \ldots, s_{i-1}, s_i a, s_{i+1}, \ldots, s_n) : a \in A\}$.
For each $s^{(a)} = (s_1, \ldots, s_{i-1}, s_i a, s_{i+1}, \ldots, s_n)$ we
have $\, |(s_1, \ldots, s_{i-1}, s_i a, s_{i+1}, \ldots, s_n)_i|$  $=$ 
$|s_i| +1 \le \ell$, and $N(s^{(a)}) = N(s) -1 < N(s)$.  
So by induction, from each
$(s_1, \ldots, s_{i-1}, s_i a, s_{i+1}, \ldots, s_n)$ one can reach
$\,\{(s_1, \ldots, s_{i-1}, s_i a, s_{i+1}, \ldots, s_n)\} \vee n A^{\ell}\,$
by one-step restrictions. Hence, 
$\, \{(s_1, \ldots, s_{i-1}, s_i a, s_{i+1}, \ldots, s_n): a \in A\}$ 
$\vee$    $n A^{\ell} \,$ can be reached from $s$ by one-step restrictions. 
Moreover, 

\smallskip

 \ \ \  \ \ \          
$\{(s_1, \ldots, s_{i-1}, s_i a, s_{i+1}, \ldots, s_n) : a \in A\}$ 
$\vee$  $n A^{\ell} \ = \ \{s\} \vee n A^{\ell}$.

\smallskip

\noindent Indeed, if $\, |s_i a| \le \ell$ and 
$\, {\rm maxlen}(s) \le \ell$ then  

\smallskip

 \ \ \  \ \ \ $s \cdot nA^* \ \cap \  n A^{\ell}$  $ \ = \ $
$\{(s_1, \ldots, s_{i-1}, s_i a, s_{i+1}, \ldots, s_n): a \in A\}$  
$\cdot$    $nA^*$   $ \ \cap \ $  $n A^{\ell}$  

 \ \ \  \ \ \ $= \ $   $\{z \in n A^{\ell} : z \ge_{\rm init} s\}$.

\smallskip

\noindent Now we apply Lemma \ref{SetJoin}(1) to express $\cap$ in terms of
$\vee$.
 \ \ \ $\Box$

\begin{lem} \label{LEMequivFin1step}
 \ Let $P, Q \subseteq n A^*$ be any finite sets, and let
$\,\ell = {\rm maxlen}(P \cup Q)$. Then the following are equivalent:

\smallskip

\noindent {\small \rm (1)} \ \ $P \equiv_{\rm fin} Q$ 
 \ \ \ {\small \rm (Def.\ \ref{DEFequivFin});}

\smallskip

\noindent {\small \rm (2)} \ \ $P \vee Q \vee n A^{\ell}$ can be 
reached 
from $P$ (and from $Q$) by finite sequences of one-step restrictions;
 
\smallskip

\noindent {\small \rm (3)} \ \ $P$ and $Q$ can be reached from one another 
by finite sequences of one-step restrictions and one-step extensions (note 
that both restriction- and extension-steps are allowed here).
\end{lem}
{\sc Proof.} $[(2) \Rightarrow (3)]$ \ This is straightforward. \\  
$[(3) \Rightarrow (1)]$ \ This follows immediately from Lemma 
\ref{OneStepRestrjoinless}(2). 

\noindent $[(1) \Rightarrow (2)]$ \ It is straightforward to prove that
${\rm maxlen}(P \cup Q) = {\rm maxlen}(P \vee Q)$.

From Lemma \ref{LEMequivStoSnAL} it 
follows that from $P$ one can reach $P \vee n A^{\ell}$; 
and from $Q$ one can reach $Q \vee n A^{\ell}$ by one-step restrictions.
By Lemma \ref{LEMequivFnAell}, $P \vee n A^{\ell} = Q \vee n A^{\ell}$, and 
hence $\,P \vee n A^{\ell}$ $=$  $P \vee Q \vee n A^{\ell}$  $=$ 
$Q \vee n A^{\ell}$. Thus, from $P$ one can reach $P \vee Q \vee n A^{\ell}$
$=$ $\,P \vee n A^{\ell}$. 
 \ \ \ $\Box$

\bigskip

\noindent Remark about Lemma \ref{LEMequivFin1step}(3):
Compare this with the converse in Lemma \ref{LEMonestepCartProd}(1).

\begin{defn} \label{DEFjonlessComplet} {\bf (completion of a joinless 
code).}
 \ For a joinless code $Q \subseteq n A$, a {\em completion} of $Q$ is any 
{\em maximal} joinless code $C \subseteq n A$ such that $Q \subseteq C$.
\end{defn}

\begin{lem} \label{LEMmaxcomplet} {\bf (existence of a completion).}
 \ For every finite joinless code $Q \subseteq n A^*$ with $|A| \ge 2$,
there exists a {\em finite} completion ${\sf C}(Q)$ of $Q$ in $n A^*$ 
such that $\,{\rm maxlen}({\sf C}(Q)) = {\rm maxlen}(Q)$.
\end{lem}
{\sc Proof.} \ Letting $\ell = {\rm maxlen}(Q)$, we pick

\medskip

\hspace{1.1in}  ${\sf C}(Q) \ = \ Q$   $ \ \cup \ $
    $\{u \in n A^{\ell} \,: \ u \not\in Q \vee n A^{\ell}\,\}\,$.

\medskip

\noindent Let us check that ${\sf C}(Q)$ has all the required properties.

It follows immediately from the formula for ${\sf C}(Q)$ that
$Q \subseteq {\sf C}(Q)$ and
$\,{\rm maxlen}({\sf C}(Q)) = {\rm maxlen}(Q)$.

The joinless code $Q$ is an essential extension of the joinless code
$Q \vee n A^{\ell}\,$; i.e.,
 \ $(Q \vee n A^{\ell}) \cdot n A^{\ell}$   $\,\subseteq\,$
$Q \cdot n A^{\ell} \ $ and
$ \ Q \equiv_{\rm fin} Q \vee n A^{\ell}\,$.
Since $\ell = {\rm maxlen}(Q)$, we have for all $u \in A^{\ell}$: 
 \ $u \in Q \vee n A^{\ell}$
 \ iff \ $u \le_{\rm init} q$ for some $q \in Q$.
The fact that ${\sf C}(Q)$ is a maximal joinless code in $n A^*$ is
now straightforward to verify directly, but is most easily seen from the
geometric representation: Every element
$u = (u_1,\ldots,u_n) \in n A^{\ell}\,$
represents a hypercube of side-length $\ell$; the hypercube represented by
$u$ is \ {\large \sf X}$_{i=1}^n [0.u_i, \ 0.u_i + \alpha^{-\ell}[ \ $
($\,\subseteq [0,1]^n$); here, $0.u_i$ denotes a rational number in
base-$\alpha$ representation, where $\alpha = |A|$.
So $n A^{\ell}$ represents the set of all hypercubes of side-length          $\ell$.
The hypercubes in $Q \vee n A^{\ell}$ tile the part of $[0,1]^n$ that is
tiled by $Q$. And the hypercubes in
$\{u \in n A^{\ell} \,:\, u \not\in Q \vee n A^{\ell} \}\,$
tile the complement of what is tiled by $Q$ in $[0,1]^n$.
Hence $\,Q$ $\,\cup\,$
$\{u \in n A^{\ell} \,:\, u \not\in Q \vee n A^{\ell}\}\,$
is a tiling of $[0,1]^n$.
 \ \ \ $\Box$

\medskip

\noindent {\bf Remark.}
The finite completion of a non-maximal finite joinless code is not unique. 

\noindent E.g., $Q = \{(a_1,a_1)\} \subseteq 2 A_2^{\,*}\,$ has the 
completions $\,{\sf C}(Q) = 2 A_2$, as well as 
$\,Q \,\cup\, \{(\e, a_0),\,(a_0,a_1)\}$,
$ \ Q \,\cup\, \{((a_0,\e),\, (a_1,a_0)\}$, and infinitely many others.

\begin{defn} \label{DEFcomplCode} {\bf (complementary joinless codes).} 
 \ Let $Q \subseteq n A^*$ be a joinless code. A {\em complementary 
joinless code} of $Q$ is any joinless code $Q' \subseteq n A^*$ such that:

\smallskip

{\small \bf (1)} \ \ $Q \cup Q'\,$ is a maximal joinless code (i.e., 
$Q \cup Q'$ is a completion of $Q$);

\smallskip

{\small \bf (2)} \ \ $Q \cdot n A^* \ \cap \ Q' \cdot n A^*$  $\,=\,$
 $\varnothing$.
\end{defn}
Note that (2) is equivalent to $\, Q \vee Q' = \varnothing$.

Def.\ \ref{DEFcomplCode} is a generalization of complementary prefix codes 
in $A^*$, defined in \cite[Def.\ 5.2]{BiCoNP} and 
\cite[Def.\ 3.29]{LRBiThompsMon}.

\begin{lem} \label{LEMcompleOmega} 
 \ Let $Q, Q' \subseteq n A^*$ be finite joinless codes. Then
$Q$ and $Q'$ are complementary joinless codes of each other \ iff 
 \ $Q \cdot n A^{\omega}$ and $Q' \cdot n A^{\omega}$ are complements as 
sets in $n A^{\omega}$.
\end{lem}
{\sc Proof.} By Lemma \ref{LEMmaxOmega}(2), $Q \cup Q'$ is a maximal joinless 
code iff $\,(Q \cup Q') \cdot nA^{\omega} = nA^{\omega}\,$.
Moreover, $Q \cdot n A^* \ \cap \ Q' \cdot n A^* \,=\, \varnothing\,$ is
equivalent to $\, Q \cdot n A^{\omega} \,\cap\, C' \cdot n A^{\omega}$ $=$
$\varnothing \ $ (by Lemma \ref{LEMdisjointOmeaga}).
 \ \ \ $\Box$

\begin{cor} \label{CORcomplCode} {\bf (existence of complementary joinless
codes).} 

\noindent Let $Q \subseteq n A^*$ be any finite joinless code, and let 
${\sf C}(Q)$ be any finite completion of $Q$.
Then $Q$ has a finite complementary joinless code $Q'$ such that
 \ $Q' = {\sf C}(Q) \minus Q\,$. \ Moreover,
$\,{\rm maxlen}(Q') = {\rm maxlen}(Q)\,$. 
\end{cor}
{\sc Proof.} This follows immediately from \ref{LEMmaxcomplet}. 
  \ \ \  $\Box$

\medskip

\noindent It is easy to see that a complementary finite joinless code of
a finite joinless code is not unique.  
 
\bigskip

\noindent Geometric meaning of complementary joinless codes: When a joinless
code $Q$ represents a set of hyperrectangles in $[0,1]^n$, $Q'$ represents a
set of hyperrectangles that tile the complement of the space tiled by $Q$.
Intuitively it is clear that a complementary joinless code code exists.

\begin{lem} \label{ComplemEquiv}
 \ Let $P, P' \subseteq n A^*$ be any finite joinless codes such that $P$ and
$P'$ are complementary.
If $Q, Q' \subseteq n A^*$ are any finite joinless codes such that
$P \equiv_{\rm fin} Q$ and $P' \equiv_{\rm fin} Q'$, then $Q$ and $Q'$
are also complementary.
\end{lem}
{\sc Proof.} Since $P \equiv_{\rm fin} Q$ and $P' \equiv_{\rm fin} Q'$, we
have by Lemma \ref{LEMequivFinOmega}(1)(2):
$\, P \cdot nA^{\omega} = Q \cdot nA^{\omega}$ and
$P' \cdot nA^{\omega} = Q' \cdot nA^{\omega}$.

Hence, $Q \cdot nA^{\omega} \,\cap\, Q' \cdot nA^{\omega}$ $=$
$P \cdot nA^{\omega} \,\cap\, P' \cdot nA^{\omega}$  $=$  $\varnothing$.
Hence, by Lemma \ref{LEMdisjointOmeaga}, we also have
$Q \cdot nA^* \,\cap\, Q' \cdot nA^* = \varnothing$.
Since $Q$ and $Q'$ are joinless, it follows that $Q \cup Q'$ is
joinless; indeed, if $q \vee q'$ existed for some $q \in Q$ and
$q' \in Q'$, then we would have $q \vee q' \in $
$Q \cdot nA^* \,\cap\, Q' \cdot nA^*$, contradicting
$Q \cdot nA^* \,\cap\, Q' \cdot nA^* = \varnothing$.

Moreover, $(Q \cup Q') \cdot nA^{\omega}$ $=$
$Q \cdot nA^{\omega} \,\cup\, Q' \cdot nA^{\omega}$ $=$
$P \cdot nA^{\omega} \,\cup\, P' \cdot nA^{\omega}$  $=$  $n A^{\omega}$.
Hence, by Lemma \ref{LEMmaxOmega},  $Q \cup Q'$ is maximal.

Thus we proved all the conditions for $Q$ and $Q'$ to be complementary.
 \ \ \ $\Box$

\begin{pro} \label{initTojoinless} {\bf (finite sets are 
{\boldmath $\equiv_{\rm fin}$}-equivalent to joinless codes).}  

\smallskip

\noindent For every finite set $S \subseteq n A^*$ there exists a
finite {\em joinless} code $C \subseteq n A^*$ such that:

\smallskip

\noindent {\small \rm (1)} \ $S \equiv_{\rm fin} C$ \ \ and
 \ \ $C \subseteq S \cdot n A^*\,$;

\smallskip

\noindent {\small \rm (2)} \ ${\rm maxlen}(C) = {\rm maxlen}(S)\,$;

\smallskip

\noindent {\small \rm (3)} \ $C$ can be reached from $S$ by a sequence of 
one-step restrictions;

\smallskip

\noindent {\small \rm (4)} \ every element of $S$ is tiled by a subset 
of $C$.
\end{pro}
Geometrically, the Proposition says that for every {\em set} of 
hyperrectangles in $[0,1]^n$, the space covered by these hyperrectangles 
can be {\em tiled} by a refinement of the given rectangles (in a tiling 
the tiles don't overlap). 
The smallest side-length of any hyperrectangle of $C$ is the same in $S$.
(3) means that $C$ is obtained from $S$ by successive $k$-partitions
along coordinate axes.

\medskip

\noindent {\sc Proof.} Letting $\, \ell = {\rm maxlen}(S)$, we choose

\smallskip

\hspace{1.in}   $C \,=\, S \vee n A^{\ell}\,$. 

\smallskip

\noindent This choice of $C$ and the proof of its properties are very 
similar to the proof of Lemma \ref{LEMmaxcomplet}.

We have $C \subseteq n A^{\ell}$. Indeed, if $s \vee u$ exists (for
$s \in S$ and $u \in n A^{\ell}$), then $s \le_{\rm init} u$, because 
$s_i \,\|_{\rm pref}\, u_i$ implies $s_i \le_{\rm pref} u_i$ when
$|u_i| = \ell = {\rm maxlen}(S)$. 

It follows from $C \subseteq n A^{\ell}$ that ${\rm maxlen}(C) = \ell$, and 
that $C$ is joinless (being a subset of the joinless code $n A^{\ell}$).

For every $s \in S$ there exists $u \in n A^{\ell}$ such that 
$s \le_{\rm init} u$ (such a $u$ is obtained from $s$ by lengthening each
$s_i$ to a string of length $\ell$); then $u = s \vee u$.  
So $s$ is the join of $\{u \in n A^{\ell} : s \le_{\rm init} u\}$, i.e.,
$s$ is tiled by these elements of $C$.

By Lemma \ref{LEMequivFPPnAell},
$S \equiv_{\rm fin} S \vee n A^{\ell}$ when $\ell = {\rm maxlen}(S)$.
By Lemma \ref{LEMequivStoSnAL}, $S \vee n A^{\ell}$ can be reached from
$S$ by a finite sequence of one-step restrictions.
 \ \ \ $\Box$

%%%%%%%%%%%%%%%%%%%%%%%%%%%%%%%%%%%%%%%%%%%%%%%%%%%%%%%%%%%%%%%%%%%%%%
%% Section
%%%%%%%%%%%%%%%%%%%%%%%%%%%%%%%%%%%%%%%%%%%%%%%%%%%%%%%%%%%%%%%%%%%%%%
\section{Right-ideal morphisms of {\boldmath $nA_k^{\,*}$} and the monoid
  {\boldmath $n M_{k,1}$} }

A {\em right-ideal morphism of $n A^*$} is any function $f$: 
$n A^* \to n A^*$ such that for all $x \in {\rm Dom}(f)$ and all 
$v \in n A^*$: \ $f(x v) = f(x) \ v$.  
(Recall that ``function'' means partial function.)
It follows that ${\rm Dom}(f)$ and ${\rm Im}(f)$ are right ideals. 
More generally, if $R$ is a right ideal then $f(R)$ and $f^{-1}(R)$ are 
right ideals; the proof is straightforward.  

Let ${\rm domC}(f)$ and ${\rm imC}(f)$ be the initial-factor codes that 
generate ${\rm Dom}(f)$ and ${\rm Im}(f)$ respectively; by
\cite[Lemma 2.7]{BinG}, these two initial-factor codes exist and are 
unique for a given right-ideal morphism $f$.  They are called the 
{\em domain code} and the {\em image code} respectively. We will use the 
notation 

\smallskip

 \ \ \  \ \ \ ${\rm maxlen}(f)$  $\,=\,$ 
${\rm maxlen} \big( {\rm domC}(f) \,\cup\, f({\rm domC}(f)) \big)\,$.

\smallskip

If $f$ and $g$ are right-ideal morphisms of $n A^*$ then the {\em composite}
$g \circ f(.)$, also denoted by $gf(.)$, is also a right-ideal morphism. 
Indeed, for all $x \in {\rm Dom}(gf)$ and $v \in nA^*$: 
$\,gf(x v) = g(f(x) \, v) = gf(x) \ v$; the latter equality holds wherever
$gf(x)$ is defined.  

We use the notation $\,f^{-1}(S) = \{x \in {\rm Dom}(f) : f(x) \in S\}\,$ 
for any set $S$ and any function $f$.

\begin{lem} \label{finvOFimCINdomC}
 \ For every right-ideal morphism $f$ of $\, n A^*:$ 

\smallskip

 \ \ \ $f^{-1}({\rm imC}(f)) \,\subseteq\, {\rm domC}(f)$;  \ \ hence 

\smallskip

 \ \ \ ${\rm imC}(f) \,\subseteq\, f({\rm domC}(f))$.

\smallskip

\noindent It follows that if $\, {\rm domC}(f)$ is finite then 
$\,{\rm imC}(f)$ is finite.
\end{lem}
{\sc Proof.} If $x \in f^{-1}({\rm imC}(f))$ then $f(x) \in {\rm imC}(f)$; 
so $f(x)$ is defined, i.e., $x \in {\rm Dom}(f)$. Therefore $x = x_0 \, v\,$
for some $x_0 \in {\rm domC}(f)$, $v \in nA^*$; hence $f(x) = f(x_0) \ v$. 
Since $f(x) \in {\rm imC}(f)$, and $f(x_0) \in {\rm Im}(f)$, and since 
${\rm imC}(f)$ is an initial-factor code, is follows that 
$v = (\e)^n$. Hence $x = x_0 v = x_0 \in {\rm domC}(f)$. 

To obtain the second inclusion we apply $f$ to 
$f^{-1}({\rm imC}(f)) \subseteq {\rm domC}(f)$. This yields
${\rm imC}(f) \subseteq f({\rm domC}(f))\,$ (since for any function $f$, 
$\,f \circ f^{-1}(.)$ is the identity function on ${\rm Im}(f)$).
 \ \ \ $\Box$

\medskip

\noindent Example where the inclusions in Lemma \ref{finvOFimCINdomC} are 
strict: \ $n=1$, $A = \{a_0, a_1\}$, ${\rm domC}(f) = A$, $f(a_0) = a_0$, 
$f(a_1) = a_0 a_0$, so ${\rm imC}(f) = \{a_0\} \ne \{a_0, a_0 a_0\}$ 
$=$  $f({\rm domC}(f))$.
And $f^{-1}({\rm imC}(f)) = f^{-1}(a_0) = \{a_0\} \ne \{a_0, a_1\}$ 
$=$  ${\rm domC}(f)$.

\begin{lem} \label{LEMequivImCfDomC}
 \ For every right-ideal morphism $f:$
 \ \ $f({\rm domC}(f)) \,\equiv_{\rm fin}\, {\rm imC}(f)\,$.
\end{lem}
{\sc Proof.} We have $f({\rm domC}(f)) \cdot nA^* \subseteq {\rm Im}(f)$
$=$  ${\rm imC}(f) \cdot nA^*$.
Also, ${\rm imC}(f) \subseteq f({\rm domC}(f))$, by Lemma
\ref{finvOFimCINdomC}.
Hence, $f({\rm domC}(f)) \cdot nA^*$  $=$  ${\rm imC}(f) \cdot nA^*$,
so $f({\rm domC}(f)) \cdot nA^* \vartriangle {\rm imC}(f) \cdot nA^*$  $=$
$\varnothing$,   which is of course finite.
 \ \ \ $\Box$

\begin{lem} \label{joinsINimC}
 \ For every right-ideal morphism $f$ of $\, n A^*:$

\smallskip

\noindent {\small \rm (1)} \ If $\,f({\rm domC}(f))$ is joinless then 
${\rm domC}(f)$ is joinless.
 \ The converse does not hold in general.

\smallskip

\noindent {\small \rm (2)} \ If $C \subseteq {\rm Im}(f)$ is an 
initial-factor code then $\,f^{-1}(C)\,$ is an initial-factor code. 

\smallskip

\noindent {\small \rm (3)} \ If $C \subseteq {\rm Im}(f)$ is a joinless code
then $\,f^{-1}(C)\,$ is joinless.

 \ The converse does not hold in general: 
There exists a right-ideal morphism $f$ of $\, n A^*$ such that 
$\,{\rm imC}(f)$ is {\em not} joinless, but $f^{-1}({\rm imC}(f))$ is 
joinless. 
Moreover, such an $f$ can be chosen so that $\,{\rm domC}(f)$ is finite and
joinless, and $f^{-1}({\rm imC}(f)) = {\rm domC}(f)$.)

\smallskip

\noindent {\small \bf (4)} \ If ${\rm domC}(f)$ is finite, then for all
$y \in {\rm Im}(f)$: \ $f^{-1}(y)$ is finite.
\end{lem}
{\sc Proof.} (1) By contradiction, if ${\rm domC}(f)$ is not joinless, 
then $a \vee b$ exists for some $a, b \in {\rm domC}(f)$ with $a \ne b$. 
Then $f(a) \vee f(b)$ also exists. Indeed, if $a \vee b = au = bv$ for some 
$u, v \in n A^*$, then $f(a) \, u = f(b) \, v$. 
Moreover, in that case, $f(a) \ne f(b)$; indeed, if we had $f(a) = f(b)$,
then $f(a) \, u = f(a) \, v$ would imply $u = v$ (since $n A^*$ is 
cancellative), and then $au = bu$ would imply $a = b$ (contradiction the 
assumption that $a \ne b$).

Counter-example for the converse: \ $n = 1$, $A = \{a_0, a_1\}$,
${\rm domC}(f) = A$, and $f(a_0) = a_0$, $f(a_1) = a_0 a_0$.

\smallskip

\noindent (2) By contradiction, if there are $a, b \in f^{-1}(C)$ with 
$a <_{\rm init} b = au$ for some $u \in n A^*$, then 
$f(a), f(b) \in C \subseteq {\rm Im}(f)$ satisfy $f(b) = f(a) \, u$. 
Since $a \ne b$, we have $u \ne (\e)^n$, hence 
$f(b) = f(a) \, u \ne f(a)$. This implies $f(a) <_{\rm init} f(b)$, 
contradicting the assumption that $C$ is an initial-factor code.  

\smallskip

\noindent (3) If there are $a, b \in f^{-1}(C)$ with $a \ne b$ such that 
$a \vee b$ exists, then $a \vee b = au = bv$ for some $u, v \in n A^*$.
Moreover, $u \ne v$ since $a \ne b$ and since $nA^*$ is a cancellative 
monoid. 
Hence $f(a), f(b) \in C \subseteq {\rm Im}(f)$ satisfy 
$f(a \vee b) = f(a) \, u = f(b) \, v$; and $f(a) \ne f(b)$ since $u \ne v$. 
This contradicts the assumption that $C$ is joinless. 

Counter-example for the converse: Let $A = \{a_0, a_1\}$, $n=2$, and
${\rm domC}(f) = \{(a_0, a_1), (a_1, a_0)\}$, which is joinless; let 
$\, {\rm imC}(f) = \{(a_0,\e), (\e, a_0)\}$, which is not 
joinless; and let $f(a_0, a_1) = (a_0,\e)$, $f(a_1, a_0) = (\e,a_0)$.
Then ${\rm imC}(f)$ is an initial-factor code that is not joinless, but
$f^{-1}({\rm imC}(f)) = {\rm domC}(f)$ is joinless. 
The right ideal $\,\{(a_0,\e), (\e, a_0)\} \cdot 2 A^*\,$ 
is not joinless generated (since the generating initial-factor code of a 
right ideal is unique, by \cite[Lemma 2.7]{BinG}). 

\smallskip

\noindent (4) For any $x \in f^{-1}(y)$ we have $f(x) = y$, and $x = pu$ 
for some $p \in {\rm domC}(f)$ and $u \in n A^*$. Hence 
$f(x) = f(p) \, u = y$. It follows that there are only finitely many choices 
for $x$, since $p \in {\rm domC}(f)$ (which is finite); and $u$ is a final
factor of $y$, and $y \in n A^*$ has only $\, \Pi_{i=1}^n (|y_i| + 1) \,$ 
final factors.   
 \ \ \ $\Box$

\bigskip

\noindent {\bf Remark.} There exists a right-ideal morphism $f$ for which 
${\rm domC}(f)$ (and hence ${\rm imC}(f)$) is finite, but {\em neither 
$\,{\rm domC}(f)\,$ nor $\,{\rm imC}(f)\,$ is joinless}. 

For example, let $n=2$, $A = \{a_0, a_1\}$, 
$\,{\rm domC}(f) = \{(a_0 a_0, a_1), (a_0, a_1a_1)\}$.  Then we have:
$\,(a_0a_0, a_1) \vee (a_0, a_1a_1) = (a_0a_0, a_1a_1)$ $=$ 
$(a_0a_0, a_1) \cdot (\e, a_1) = $ $(a_0, a_1a_1) \cdot (a_0, \e)$; so 
${\rm domC}(f)$ is not joinless. 
And $f$ satisfies $\, f(a_0a_0, a_1a_1) = f(a_0a_0, a_1) \cdot (\e, a_1) =$
$f(a_0, a_1a_1) \cdot (a_0, \e)$. This puts limitations on how 
$f(a_0 a_0, a_1)$ and $f(a_0, a_1a_1)$ can be chosen (for the given 
${\rm domC}(f)$).  E.g., we can choose
$\,f(a_0 a_0, a_1) = (x a_0, y)$ and $f(a_0, a_1a_1)  = (x,y a_1)$ 
for any $x, y \in A^*$. In that case, 
${\rm imC}(f) = \{(x a_0,y), (x,y a_1)\}$; 
and $(x a_0,y) \vee (x,y a_1) = (x a_0, y a_1)$ exists.

\bigskip

\noindent In \cite{BinG} the string-based definition of the Brin-Thompson 
group used the  inverse monoid 

\medskip

$n {\cal RI}^{\sf fin}$  $ \ = \ $
$\{f : \, f$ is a right ideal morphism of $nA^*$ such that $f$ is injective,

\hspace{1.06in} and ${\rm domC}(f)$ and ${\rm imC}(f)$
are {\em finite, maximal, joinless} codes\} .

\smallskip

\noindent There are a number of ways to generalize this:

%%%%%%%%%%% n RM
\begin{defn} \label{defnRMfin}  {\bf (sets and monoids of right-ideal 
morphisms.)} 

\smallskip

\noindent For any $n \ge 1$ and any alphabet $A$, we define the following 
sets of right-ideal morphisms of $n A^*$: 

\medskip
 
$n\,{\cal RM}^{\sf fin}_0 \ = \ \{f : f$ is a right-ideal morphism
of $n A^*$ such that ${\rm domC}(f)$ is {\em finite}\};

\medskip

$n\,{\cal RM}^{\sf fin}_1 \ = \ \{f : f$ is a right-ideal morphism
of $n A^*$ such that ${\rm domC}(f)$ is {\em finite} and {\em joinless}\};

\medskip

$n\,{\cal RM}^{\sf fin}_2 \ = \ \{f : f$ is a right-ideal morphism
of $n A^*$ such that both ${\rm domC}(f)$ and ${\rm imC}(f)$ are 

\hspace{1.15in} {\em finite} and {\em joinless}\};

\medskip

$n\,{\cal RM}^{\sf fin}_{\rm norm} \ = \ \{f $ 
$\in n\,{\cal RM}^{\sf fin}_2 : \ f({\rm domC}(f)) = {\rm imC}(f)\}$.
\end{defn}

\medskip

\noindent Obviously, \ $n\,{\cal RM}^{\sf fin}_{\rm norm} \subseteq$
$n\,{\cal RM}^{\sf fin}_2 \subseteq n\,{\cal RM}^{\sf fin}_1$
$\subseteq n\,{\cal RM}^{\sf fin}_0$; and it is not hard to prove that the 
inclusions are strict.
The elements of $n\,{\cal RM}^{\sf fin}_{\rm norm}$ are called {\em normal}
(see \cite[Def.\ 5.6]{BiThompsMonV3}, \cite[Def.\ 4.5A]{Equiv}). 

For for every $f \in n\,{\cal RM}^{\sf fin}_0$: $ \ f({\rm domC}(f))$ is
finite; hence by Lemma \ref{finvOFimCINdomC}, ${\rm imC}(f)$ is finite.

\medskip

Every $f \in n\,{\cal RM}^{\sf fin}_0$ can be described by a unique 
finite {\bf table}; this is the restriction $f|_P: P \twoheadrightarrow Q$ 
of $f$, where $P = {\rm domC}(f)$, and $Q = f({\rm domC}(f))$. Here
$P$ can be any finite initial-factor code (which is joinless in case 
$f \in  n\,{\cal RM}^{\sf fin}_1$). And $Q$ ($= f(P)$) can be any finite 
subset of $n A^*$ of size $\le |P|$.
Moreover, the right ideal $\,Q \cdot n A^*\,$ is generated by the finite 
initial-factor code ${\rm imC}(f)$. By Lemma \ref{finvOFimCINdomC}, 
$\,{\rm imC}(f) \subseteq Q$; in general, ${\rm imC}(f)$ can be a strict 
subset of $Q$.

\begin{pro} \label{TableCompos}
 \ Let $f$ and $g$ be right-ideal morphisms of $n A^*$. \ Then: 

\smallskip

\noindent {\small \bf (1)} \ \ ${\rm domC}(gf(.)) \ = \ $
${\rm domC}(f) \,\vee\, f^{-1}({\rm domC}(g))$.

\smallskip

\noindent {\small \bf (2)} \ \ If ${\rm domC}(f)$ and ${\rm domC}(g)$ 
are joinless codes then $\,{\rm domC}(gf(.))\,$ is joinless. 
\end{pro}
{\sc Proof.} (1) We have: $\,x \in {\rm Dom}(gf(.))$ iff $\,gf(x)$ is 
defined. 
The latter holds iff $x \in {\rm Dom}(f)$ and $f(x) \in {\rm Dom}(g)$;
hence, ${\rm Dom}(gf(.)) = {\rm Dom}(f) \,\cap\, f^{-1}({\rm Dom}(g))$.
The generating initial factor code of ${\rm Dom}(f)$ is ${\rm domC}(f)$; and
by Lemma \ref{joinsINimC}(2), the generating code of $f^{-1}({\rm Dom}(g))$ 
is $\,f^{-1}({\rm domC}(g))$.
Therefore, by Lemma \ref{SetJoin} the generating code of 
$\,{\rm Dom}(f) \,\cap\, f^{-1}({\rm Dom}(g))\,$ is 
$\, {\rm domC}(f) \,\vee\, f^{-1}({\rm domC}(g))$.
Since the generating initial factor code of a right ideal is unique (by 
\cite[Lemma 2.7]{BinG}), we conclude that
$\,{\rm domC}(gf(.)) \,=\, {\rm domC}(f) \,\vee\, f^{-1}({\rm domC}(g))$.

(2) By Lemma \ref{joinsINimC}(3), $f^{-1}({\rm domC}(g))$ is joinless if
${\rm domC}(g)$ is joinless. By Lemma \ref{SetJoin} the result then follows.
 \ \ \ $\Box$

\begin{pro} \label{nRNfinMon} {\bf (monoids).} 
 \ The sets $\,n\,{\cal RM}^{\sf fin}_0$ and $\,n\,{\cal RM}^{\sf fin}_1$ 
are monoids under composition.
\end{pro}
{\sc Proof.} We saw that if $f,g$ are right-ideal morphisms then $fg(.)$ is 
a right-ideal morphism.
If ${\rm domC}(f)$ and ${\rm domC}(g)$ are finite, then ${\rm domC}(gf)$ is
finite, by Prop.\ \ref{TableCompos}(1) and Lemma \ref{joinsINimC}(2). 
Hence $\,n\,{\cal RM}^{\sf fin}_0$ is a monoid.

If ${\rm domC}(f)$ and ${\rm domC}(g)$ are joinless codes then 
${\rm domC}(gf)$ is also joinless, by Prop.\ \ref{TableCompos}(2).
Hence $\,n\,{\cal RM}^{\sf fin}_1$ is a monoid.
 \ \ \ $\Box$

\begin{lem} \label{nRNfin_2notMOn}
 \ Let $n \ge 2$ and $A = \{a_0, a_1\}$. 

\smallskip

\noindent {\small \rm (1)} \ There exist $f \in n\,{\cal RM}^{\sf fin}_2$ 
and a finite joinless code $C \subseteq {\rm Dom}(f) \subseteq n A^*$ such 
that $\, f(C) \cdot n A^*\,$ is not joinless generated.

\smallskip

\noindent {\small \rm (2)} \ The set $\,n\,{\cal RM}^{\sf fin}_2$ is 
{\em not} closed under composition.
\end{lem}
{\sc Proof.} (1) Let $f =$ 
$\{\big((a_0, a_0), (\e, \e)\big), \,$  
$\big((a_0, a_1),(\e, a_0)\big),\,$  
$\big((a_1, a_0), (a_0,\e)\big)\}$. 
So $\,{\rm domC}(f) = $ $\{(a_0,a_0),$  $(a_0,a_1),$  $(a_1,a_0)\}$, which 
is joinless.  And $f({\rm domC}(f)) = $
$\{(\e, \e),$  $(\e, a_0),$  $(a_0,\e)\}$; hence, 
${\rm imC}(f) = \{(\e, \e)\}$, which is joinless.

Let $C = \{(a_0, a_1),(a_1,a_0)\}$, which is joinless. Then $f(C) =$ 
$\{(\e,a_0),(a_0,\e)\}$, which is an initial-factor code that 
is not joinless. Moreover, the right ideal 
$\,\{(\e,a_0),(a_0,\e)\} \cdot 2 A^*$ is not joinless 
generated, since it is generated by $\{(\e,a_0),(a_0,\e)\}$,
which (by Lemma \cite[Lemma 2.7]{BinG}) is the unique initial-factor code 
that generates this right ideal.

(2) Let $f_1 = {\sf id}_{C \cdot 2 A^*}$, i.e., the restriction of the 
identity function, where $C$ is the joinless code above. And let
$f_2 = f$, as defined in (1). Then $f_1, f_2 \in 2\,{\cal RM}^{\sf fin}_2$.  
And ${\rm domC}(f_2 f_1(.)) = C$, which is joinless,
but $f({\rm domC}(f_2 f_1(.))) = {\rm imC}(f_2 f_1(.)) = f(C)$ $=$ 
$\{(\e,a_0),(a_0,\e)\}$, which is not joinless. 
 \ \ \ $\Box$

\begin{lem} \label{nRNfin_NormnotMon}
 \ For all $n \ge 1$ and $A$ with $|A| \ge 2$, the set 
$\,n\,{\cal RM}^{\sf fin}_{\rm norm}$ is {\em not} closed under composition.
\end{lem}
{\sc Proof.} Here is an example where $n = 1$ and $A = \{a_0, a_1\}$, taken 
from \cite[Prop.\ 5.8]{Equiv}; this example can easily be extended to 
examples with any $n > 1$ and any $A$ with $|A| > 2$.  

Consider $f,g \in {\cal RM}^{\sf fin}$ with 
${\rm domC}(f) = A = {\rm domC}(g)$, and $f(a_0) = a_0, \ f(a_1) = a_1 a_0$,
 \ $g(a_0) = g(a_1) = a_0$. Then 
$f({\rm domC}(f)) = {\rm imC}(f) = \{a_0, a_1 a_0\}$, 
and $g({\rm domC}(g)) = {\rm imC}(g) = \{a_0\}$, hence $f$ and $g$ are normal.
But ${\rm domC}(gf) = A$, $gf(a_0) = a_0, \ gf(a_1) = a_0 a_0$, thus
$gf({\rm domC}(gf)) = \{a_0, a_0a_0\}$. So $\,gf({\rm domC}(gf))$ is not a 
prefix code, hence $gf$ is not normal. 
 \ \ \ $\Box$

\bigskip

In \cite[Lemma 5.1]{Equiv} it was proved that if $f \in {\cal RM}^{\sf fin}$
is injective then it is normal; in the next Lemma we prove the same thing for
$n {\cal RM}_0^{\sf fin}$. Therefore, non-normal morphisms do not appear in the
study of the Brin-Thompson groups $n V$.

\begin{lem} \label{injNormal}
 \ If $g$ is an injective right-ideal morphism of $n A^*$,
then $g$ is normal.
\end{lem}
{\sc Proof.} Let $x \in {\rm domC}(g)$.  Then $g(x) \in {\rm Im}(g)$, so 
$g(x) = uv$ for some $u \in {\rm imC}(g)$, $v \in n A^*$. 
Let $z \in {\rm Dom}(g)$ be such that $g(z) = u$; then
$z = st$ for some $s \in {\rm domC}(g)$, $t \in n A^*$.
Hence, $g(x) = uv = g(z) \, v = g(zv) = g(stv)$. Since $g$ is injective, this
implies that $x = stv$, hence $s \le_{\rm init} x$. But then $s = x$, since
$x$ and $s$ belong to the initial-factor code ${\rm domC}(g)$. Therefore
$t = v = (\e)^n$. It follows that $g(x) = uv = u \in {\rm imC}(g)$,
so $g(x) \in {\rm imC}(g)$. Thus, $g({\rm domC}(g)) \subseteq {\rm imC}(g)$.
  \ \ \ $\Box$

\begin{defn} \label{DEFequivMorph} {\bf ({\boldmath $\,\equiv_{\rm fin}$} 
between morphisms).} 
 \ The relation $\equiv_{\rm fin}$ between right-ideal morphisms 
$f$ and $g$ of $\,n A^*$ is defined as follows:

\smallskip

 \ \ \ $f \equiv_{\rm fin} g$ \ \ \ iff
 \ \ \ ${\rm domC}(f) \equiv_{\rm fin} {\rm domC}(g)$ \ and
 \ $f|_{{\rm Dom}(f) \,\cap\, {\rm Dom}(g)}$  $=$ 
$g|_{{\rm Dom}(f) \,\cap\, {\rm Dom}(g)}$. 
\end{defn}

\begin{defn} \label{DEFmorphOmega} 
 \ Let $f$ be any right-ideal morphism of $\,n A^*$. Then $f$ can be 
extended to $A^{\omega}$ as follows:

\medskip

 \ \ \ ${\rm Dom}(f|_{n A^{\omega}})$  $\,=\,$ 
  ${\rm domC}(f) \cdot n A^{\omega}\,$, \ \ \ and
 
\medskip

  \ \ \ $f|_{n A^{\omega}}(p u) = f(p) \ u$, \ for 
  $\,p \in {\rm domC}(f)$ and $\,u \in n A^{\omega}$.
\end{defn}
Henceforth, $f$ continues to denote the application of $f$ to $\,n A^*$;
the application of $f$ to $n A^{\omega}$ is denoted by $f|_{n A^{\omega}}$.

\begin{lem} \label{LEMequivMorphOmega} {\bf ({\boldmath $\,\equiv_{\rm fin}$
and $nA^{\omega}$}).}
 \ For any right-ideal morphisms $f$ and $g$ of $\,n A^*$ with finite 
domain codes we have:

\smallskip

\hspace{0.8in}  $f \equiv_{\rm fin} g$ \ \ iff 
 \ \ $f|_{nA^{\omega}} = g|_{nA^{\omega}}\,$.
\end{lem}
{\sc Proof.} Recall that by Lemma \ref{LEMequivFinOmega}, 
${\rm domC}(f) \equiv_{\rm fin} {\rm domC}(g)\,$ iff 
$\,{\rm domC}(f) \cdot n A^{\omega} = {\rm domC}(g) \cdot n A^{\omega}$.
By Def.\ \ref{DEFmorphOmega}, this is also equivalent to 
$ \ {\rm Dom}(f|_{nA^{\omega}})$  $=$ ${\rm Dom}(g|_{nA^{\omega}})$.
 
\smallskip

\noindent $[\Rightarrow]$ If $f \equiv_{\rm fin} g$ then, as we just saw,
$\,{\rm Dom}(f|_{nA^{\omega}})$  $=$ ${\rm Dom}(g|_{nA^{\omega}})$.
Hence for every $w \in {\rm Dom}(f|_{nA^{\omega}})$ we have: 
$w \in {\rm domC}(f) \cdot nA^{\omega}$  $\,\cap\,$ 
${\rm domC}(g) \cdot nA^{\omega}$ $=$  
$({\rm domC}(f) \vee {\rm domC}(g)) \cdot nA^{\omega}$; so, $w = z u$ for 
some $z \in {\rm domC}(f) \vee {\rm domC}(g)$ and $u \in nA^{\omega}$. 
Since, by definition, $f \equiv_{\rm fin} g$ implies $f(x) = g(x)$ for all
$x \in {\rm Dom}(f) \,\cap\, {\rm Dom}(g)$, we conclude: 
$f(z) = g(z)$. Hence, $f(w) = g(w)$.
Hence, $f|_{nA^{\omega}}(w) = g|_{nA^{\omega}}(w)$.

\smallskip

\noindent $[\Leftarrow]$ If $f|_{nA^{\omega}} = g|_{nA^{\omega}}$ then
they have the same domain, so by the above, 
${\rm domC}(f) \equiv_{\rm fin} {\rm domC}(g)$. 
We still have to show that for every 
$x \in {\rm Dom}(f) \,\cap\, {\rm Dom}(g)$: $\, f(x) = g(x)$.

Since $f|_{nA^{\omega}} = g|_{nA^{\omega}}$ we have for every 
$x \in {\rm Dom}(f) \,\cap\, {\rm Dom}(g)$ and every
$u \in nA^{\omega}$: $\, f|_{nA^{\omega}}(x u) = g|_{nA^{\omega}}(x u)$.
Hence, by the definition of $f|_{nA^{\omega}}$ and $g|_{nA^{\omega}}:$
$\, f|_{nA^{\omega}}(x u) = f(x) \, u = g(x) \, u = g|_{nA^{\omega}}(x u)$.
For a given $x$, this holds for all $u \in nA^{\omega}$; hence, 
$f(x) = g(x) \ $ (by Lemma \ref{LEMprefEquOmega}).
 \ \ \ $\Box$

\begin{lem} \label{LEMequivMorphTrans}
 \ The relation $\equiv_{\rm fin}$ on $\, n \,{\cal RM}^{\sf fin}_0$ 
is transitive.
\end{lem}
{\sc Proof.}  This follows immediately from Lemma
\ref{LEMequivMorphOmega}.
 \ \ \ $\Box$

\begin{lem} \label{LEMequivOimC}
 \ For any right-ideal morphisms $f, g$ of $nA^*$ with finite domain codes,
 \ $f \equiv_{\rm fin} g$  \ implies
 \ ${\rm imC}(f) \equiv_{\rm fin} {\rm imC}(g)$.
\end{lem}
{\sc Proof.} By Lemma \ref{LEMequivMorphOmega}, $f \equiv_{\rm fin} g$ 
implies ${\rm domC}(f) \cdot nA^{\omega} = {\rm domC}(g) \cdot nA^{\omega}$.
Hence, $f({\rm domC}(f) \cdot nA^{\omega})$  $=$ 
$g({\rm domC}(g) \cdot nA^{\omega})$. 
Moreover, $f({\rm domC}(f) \cdot nA^{\omega}) $  $=$ 
$f({\rm domC}(f)) \cdot nA^{\omega} \ $ (and similarly for $g$), since 
$f$ is a right-ideal morphism, and ${\rm domC}(f) \subseteq {\rm Dom}(f)$.
Hence, $f({\rm domC}(f)) \cdot nA^{\omega}$  $=$ 
$g({\rm domC}(g)) \cdot nA^{\omega}$, so 
$f({\rm domC}(f)) \equiv_{\rm fin} g({\rm domC}(g))$.

Hence by Lemma \ref{LEMequivImCfDomC}, 
${\rm imC}(f) \equiv_{\rm fin} f({\rm domC}(f))$.
 \ \ \ $\Box$

\begin{lem} \label{LEMomegaProduct}
 \ For all right-ideal morphisms $k, h$ on $nA^*:$ 
 \ \ $(k \circ h)|_{nA^{\omega}}$  $\,=\,$ 
$k|_{nA^{\omega}} \circ h|_{nA^{\omega}}$.
\end{lem}
{\sc Proof.} By Def.\ \ref{DEFmorphOmega}, 
$ \ {\rm Dom}\big((k \, h)|_{nA^{\omega}}\big)$  $\,=\,$ 
${\rm Dom}\big(k|_{nA^{\omega}} \ h|_{nA^{\omega}}\big)$  $\,=\,$
${\rm domC}(k \, h) \cdot nA^{\omega}$.  
And for all $x = p u\,$ with $p \in {\rm domC}(k\,h)$ and 
$u \in nA^{\omega}$: 
 \ $(k \, h)|_{nA^{\omega}}(x)$  $\,=\,$  $k(h(p)) \ u$ $\,=\,$ 
$k|_{nA^{\omega}} \big( h(p) \, u \big)$  $\,=\,$
$k|_{nA^{\omega}} \big( h|_{nA^{\omega}}(p u) \big)$. 
 \ \ \ $\Box$

\begin{lem} \label{LEMequivMorphCongr}
 \ The relation $\equiv_{\rm fin}$ on $\, n \,{\cal RM}^{\sf fin}_0$
is a {\em congruence}.
\end{lem}
{\sc Proof.} The relation $\equiv_{\rm fin}$ is obviously reflexive and
symmetric, and by Lemma \ref{LEMequivMorphTrans} it is transitive.
By Lemma \ref{LEMequivMorphOmega}, the congruence property is equivalent to
the following sentence: 

\smallskip

$(\forall f, g, h \in n \,{\cal RM}^{\sf fin}_0)$
$[\, f|_{nA^{\omega}} = g|_{nA^{\omega}} \ $ implies
$ \ (hf)|_{nA^{\omega}} = (hg)|_{nA^{\omega}} \,$ and
$\, (fh)|_{nA^{\omega}} = (gh)|_{nA^{\omega}} \,]$.
 
\smallskip

\noindent This, in turn, follows from Lemma \ref{LEMomegaProduct}. 
 \ \ \ $\Box$

\bigskip

%%% DEF of Brin-Thompson group monoid 
\noindent We can now define the {\em monoid version $n M_{k,1}$ of the 
Brin-Higman-Thompson group} $n G_{k,1}$:

\begin{defn} \label{DEFnM21} {\bf (definition of {\boldmath
$\, n M_{k,1}$).}}

\smallskip

The monoid $\, n M_{k,1}\,$ is the transformation monoid of the action of 
$\,n \,{\cal RM}^{\sf fin}_1$ on $\,nA_k^{\,\omega}$, where 
$A_k = \{a_0,a_1, \ldots, a_{k-1}\}$. 
In other words, $ \ n M_{k,1}\,$  $=$ 
$\, \{f|_{nA_k^{\,\omega}} \,:\, f \in n \,{\cal RM}^{\sf fin}_1\}$.

\smallskip

Hence $\, n M_{k,1}\,$  is isomorphic to
$\,n \,{\cal RM}^{\sf fin}_1 \!/\! \equiv_{\rm fin}\,$. 
\end{defn}
The latter is the quotient monoid of the monoid $n {\cal RM}^{\sf fin}_1$ 
by the congruence $\equiv_{\rm fin}$. 
The definition of $\,nM_{k,1}\,$ uses $n {\cal RM}^{\sf fin}_1$, consisting 
of the right-ideal morphisms whose domain codes are finite and joinless, 
but in Lemma \ref{LEMhhINVdomC} we will see that every 
$f \in n {\cal RM}^{\sf fin}_0$ is $\,\equiv_{\rm fin}$-equivalent to some
$\varphi \in n {\cal RM}^{\sf fin}_{\rm norm}$; hence, $n M_{k,1}$ is the 
same if $n {\cal RM}^{\sf fin}_0$, or if $n {\cal RM}^{\sf fin}_2$, or 
$n {\cal RM}^{\sf fin}_{\rm norm}$ is used instead of 
$n {\cal RM}^{\sf fin}_1$; see Corollary \ref{nM21RM0RM}.

\begin{lem} \label{LEMequivMorphRestr}
 \ For every $f \in n \,{\cal RM}^{\sf fin}_0$ and every finite set
$P \subseteq {\rm Dom}(f)$: 

\smallskip

 \ \ \  \ \ \ $\,P \,\equiv_{\rm fin}\, {\rm domC}(f)$ \ implies
 \ $f \equiv_{\rm fin} f|_{P \, nA^*}\,$.
\end{lem}
{\sc Proof.} Obviously, $f$ and $f|_{P \, nA^*}$ agree on 
$\,{\rm Dom}(f) \cap {\rm Dom}(f|_{P \, nA^*})$ $=$ $P \, nA^*$.
This, together with $\,P \equiv_{\rm fin} {\rm domC}(f)$, implies
the Lemma.
 \ \ \ $\Box$

\begin{lem} \label{LEMequivMorphBound}
 \ For all $f, g \in n \,{\cal RM}^{\sf fin}_0$ and 
$\,\ell$  $\,\ge\,$  $\max\{{\rm maxlen}(f), {\rm maxlen}(g)\}:$

\smallskip

 \ \ \ $f \equiv_{\rm fin} g$ \ \ iff
 \ \ $f|_{n A^{\ell}} \,=\, g|_{n A^{\ell}}$ \ \ iff
 \ \ $f|_{n A^{\ge \ell}} \,=\, g|_{n A^{\ge \ell}}$ .
\end{lem}
{\sc Proof.} The second equivalence holds because $\,n A^{\ell} \subseteq$
${\rm Dom}(f) \,\cap\, {\rm Dom}(g)$, and $f$ and $g$ are right-ideal 
morphisms. Let us prove the first equivalence. 

We have: \ ${\rm domC}(f|_{n A^{\ell}})$  $\,=\,$
${\rm domC}(f) \,\vee\, n A^{\ell}$ $ \ \equiv_{\rm fin} \ $ 
${\rm domC}(f) \ $ (the latter by Lemma \ref{LEMequivFnAell}).
Moreover, $n A^{\ell} \subseteq {\rm Dom}(f) \,\cap\, {\rm Dom}(g)$; hence 
by Lemma \ref{LEMequivMorphRestr}, $f \equiv_{\rm fin} f|_{n A^{\ell}}$.
Thus, $f|_{n A^{\ell}} \,=\, g|_{n A^{\ell}}\,$ implies 
$\, f \equiv_{\rm fin} g$.

Conversely, $f \equiv_{\rm fin} g$ implies that $f$ and $g$ agree on 
${\rm Dom}(f) \,\cap\, {\rm Dom}(g)$. This implies that that $f$ and $g$ 
agree on $n A^{\ell}$, since $n A^{\ell} \subseteq$ 
${\rm Dom}(f) \,\cap\, {\rm Dom}(g)$.
 \ \ \ $\Box$

\begin{lem} \label{LEMfOFeqiv} 
 \ For all $f \in n\,{\cal RM}^{\sf fin}_0$, and all finite sets 
$P, Q \subseteq {\rm Dom}(f):$

\smallskip

 \ \ \  \ \ \ $P \equiv_{\rm fin} Q$ \ implies 
 \ $f(P) \equiv_{\rm fin} f(Q)$.
\end{lem}
{\sc Proof.} Since $P, Q \subseteq {\rm Dom}(f)$, we have 
$f(P \ nA^{\omega}) = f(P) \ nA^{\omega}$ and
$f(Q \ nA^{\omega}) = f(Q) \ nA^{\omega}$. By Lemma \ref{LEMequivFinOmega}, 
$P \equiv_{\rm fin} Q$ implies $P \ nA^{\omega} = Q \ nA^{\omega}$,
hence $\,f(P) \ nA^{\omega} = f(Q) \ nA^{\omega}$. Hence 
$f(P) \equiv_{\rm fin} f(Q)$ (by Lemma \ref{LEMequivFinOmega}).
 \ \ \ $\Box$

\begin{lem} \label{LEMinvimage}
 \ Let $h$ be any right-ideal homomorphism of $nA^*$ with finite 
${\rm domC}(h)$. Let $Q \subseteq {\rm Im}(h)$ be any finite set 
satisfying the following condition.

\medskip

\noindent $(\star)$ \  \ \ \  For all $q \in Q$ and all 
$y \in h({\rm domC}(h))$: 
  \ if $\,y \vee q\,$ exists, then $\,y \le_{\rm init} q$.  

\medskip

\noindent Then:

% \smallskip

\hspace{0.8in} 
$h^{-1}(Q \cdot nA^{\omega}) \,=\, h^{-1}(Q) \cdot nA^{\omega}$.
\end{lem}
Since $\,y \le_{\rm init} q\,$ iff $\,y \vee q = q$, condition $(\star)$
is equivalent to: \ If $\,y \vee q\,$ exists, then $\,y \vee q = q$.

Informally, condition $(\star)$ says that every element of $Q$ is
``longer'' than all join-comparable elements in the finite set 
$h({\rm domC}(h))$.

\medskip

\noindent 
{\sc Proof.} We have $\,h^{-1}(Q \cdot nA^{\omega})$  $\subseteq$ 
${\rm Dom}(h|_{nA^{\omega}})$.    Moreover, every element of 
${\rm Dom}(h|_{nA^{\omega}})$ is of the form $x u$, with
$x \in {\rm domC}(h)$ and $u \in nA^{\omega}$.

\smallskip

\noindent $[\supseteq ]$  Consider $pv \in h^{-1}(Q) \cdot nA^{\omega}$, for
any $p \in h^{-1}(Q)$ and $v \in nA^{\omega}$. Then $\,h(p v) = h(p) \, v$,
since $p \in h^{-1}(Q)$ $\subset$  ${\rm Dom}(h)$.  
And $h(p) \, v$ $\in$  $h(h^{-1}(Q)) \cdot nA^{\omega}$ $=$ 
$Q \cdot nA^{\omega}$.
Hence, $h(p v) \in Q \cdot nA^{\omega}$, which is equivalent to 
$p v \in h^{-1}(Q \cdot nA^{\omega})$.

\smallskip  

\noindent $[\subseteq]$ Consider $x u \in h^{-1}(Q \cdot nA^{\omega})$
$\subseteq$  ${\rm Dom}(h|_{nA^{\omega}})$, where $x \in {\rm domC}(h)$ 
and $u \in nA^{\omega}$; then $\,h(xu) \in Q \cdot nA^{\omega}$.
Since $x \in {\rm Dom}(h)$: $\,h(xu) = h(x) \, u$. 
And since $\, h(x) \, u \in Q \cdot nA^{\omega}$: $\, h(x) \, u = q v\,$ 
for some $q \in Q$ and $v \in nA^{\omega}$. By Lemma \ref{LEMjoinOmega} 
the latter equality implies that $\,h(x) \vee q\,$ exists.
Condition $(\star)$ in the Lemma then implies: 
 \ $h(x) \le_{\rm init} q = h(x) \, z$, for some $z \in nA^*$.  
 
Hence, $h(x) \, u = q v = h(x) \, z \, v$.  By left-cancellativity in
$nA^{\omega}$ this implies $\, u = z v$.
Therefore, $x u = x z v$; moreover, $x z v \in h^{-1}(q) \ v$, since we 
saw that $q = h(x) \, z = h(xz) \,$ (and $q = h(xz)$ is equivalent to 
$x z \in h^{-1}(q)$).  Thus, 
$x u = x z v \in h^{-1}(q) \ v \,\subseteq\, h^{-1}(Q) \cdot n A^{\omega}$. 
  \ \ \ $\Box$

\begin{lem} \label{LEMequivInvImage}
 \ Let $h$ be a right-ideal homomorphism of $nA^*$ with finite
${\rm domC}(h)$. Let $P, Q \subseteq {\rm Im}(h)$ be any finite sets 
satisfying the following condition.

\medskip

\noindent $(\star)$ \  \ \ \  For all $r \in P \cup Q$ and all
$y \in h({\rm domC}(h))$:
  \ if $\,y \vee r\,$ exists, then $\,y \le_{\rm init} r$.

\medskip

\noindent Then:

% \smallskip

\hspace{0.8in}  $P \equiv_{\rm fin} Q$ \ \ implies
 \ \ $h^{-1}(P) \,\equiv_{\rm fin}\,  h^{-1}(Q)$. 
\end{lem}
{\sc Proof.} By Lemma \ref{LEMequivFinOmega}, the conclusion is 
equivalent to: \ $P \cdot nA^{\omega}  = Q \cdot nA^{\omega} \ $ implies
$ \ f^{-1}(P) \cdot nA^{\omega}  = f^{-1}(Q) \cdot nA^{\omega}$. 
 \ Obviously, $P \cdot nA^{\omega}  = Q \cdot nA^{\omega}\,$ implies
$\, h^{-1}(P \cdot nA^{\omega}) = h^{-1}(Q \cdot nA^{\omega})$.  
By Lemma \ref{LEMinvimage}, 
$ \ h^{-1}(P) \cdot nA^{\omega} = h^{-1}\big(P \cdot nA^{\omega}\big)$, and
$ \ h^{-1}(Q) \cdot nA^{\omega} = h^{-1}(Q \cdot nA^{\omega})$.
The result then follows immediately.
 \ \ \ $\Box$

\bigskip

\noindent {\bf Remark:} There exist $f \in n\,{\cal RM}^{\sf fin}_0\,$ and 
finite sets $P, Q \subseteq {\rm Im}(f)$ such that
$ \ P \equiv_{\rm fin} Q$, but $ \ f^{-1}(P) \not\equiv_{\rm fin} f^{-1}(Q)$.
Hence, in Lemmas \ref{LEMequivInvImage} and \ref{LEMinvimage}, condition $(\star)$ (or some other non-trivial assumption) is needed.

\smallskip

\noindent Example: Let $n=1$, $A = \{a_0, a_1\}$, 
$\,{\rm domC}(f) = A$, and $\,f(a_0) = a_0$, $f(a_1) = a_0 a_0$.  
So $\, {\rm imC}(f) = \{a_0\}$  $\equiv_{\rm fin}$  
$\{a_0, a_0a_0\} = f({\rm domC}(f))$.
Moreover,  
$\,f^{-1}({\rm imC}(f)) = f^{-1}(a_0) = \{a_0\}$  $\,\not\equiv_{\rm fin}\,$
$A = {\rm domC}(f) = f^{-1}(f({\rm domC}(f)))$ 
 \ (in this example).
So, letting $P = {\rm imC}(f)\,$ and $Q = f({\rm domC}(f))$, we have
$P \equiv_{\rm fin} Q\,$ but 
$\,f^{-1}(P) \,\not\equiv_{\rm fin}\, f^{-1}(Q)$.
 \ \ \ $\Box$
 
\bigskip

\begin{defn} \label{DEFonestepMorph} {\bf (one-step restriction or 
extension of a morphism).} 
 \ Let $f \in n \,{\cal RM}^{\sf fin}_1$ be any right-ideal morphism 
with $P = {\rm domC}(f)$ a finite joinless code, let $p \in P$, let 
$i \in \{1,\ldots,n\}$, and let $P_{p,i}$ be the one-step restriction of 
$P$ at $p$ and $i$ (as in Lemma {\rm \ref{OneStepRestrjoinless}}). 

Then the restriction $f|_{P_{p,i}}$ is called a {\em one-step
restriction} of $f$; and $f$ is called a {\em one-step extension} of 
$f|_{P_{p,i}}$.
\end{defn}
A table for $f|_{P_{p,i}}$ is obtained from the table 
$f|_P: P \twoheadrightarrow Q\,$ for $f$ by replacing the entry $(p, f(p))$ 
by the set of entries 
%$\,\{\big((p_1, \,\ldots, p_{i-1},\, p_i a,\, p_{i+1},\, \ldots, p_n), $
%$ f(p_1, \,\ldots, p_{i-1},\, p_i a,\, p_{i+1},\, \ldots, p_n)\big) \, :\,$
%$a \in A\}$ $=$ 
 \ $(p, f(p))$ $\cdot$
$\big( \{\e\}^{i-1} \times A \times \{\e\}^{n-i} \big)$.

\begin{lem} \label{nRMeqiv} {\bf (closure under one-step restriction).} 
 \ If $\,f \in n\,{\cal RM}^{\sf fin}_1\,$  and $f|_R$ is reached from $f$ 
by a one-step restriction, then $f|_R \in n {\cal RM}^{\sf fin}_1$.
Similarly, $n\,{\cal RM}^{\sf fin}_2$ and 
$n\,{\cal RM}^{\sf fin}_{\rm norm}$ are closed under one-step restriction. 
\end{lem}
{\sc Proof.} If $f \in n\,{\cal RM}^{\sf fin}_1$: 
In a one-step restriction, the cardinality of the domain code 
is increased by $|A| - 1$; so finiteness of the domain code is preserved.
And by Lemma \ref{OneStepRestrjoinless}, joinlessness of the domain code
is preserved too. 

If $f \in n\,{\cal RM}^{\sf fin}_2$:
Then in the table $f|_P: P \twoheadrightarrow Q\,$ for $f$, ${\rm imC}(f)$ 
($\subseteq Q$) is joinless. 
Hence, either ${\rm imC}(f) = {\rm imC}(f|_{P_{p,i}})$, or in 
${\rm imC}(f)$, some entry $q \in Q$ is replaced by $\, q$ $\cdot$ 
$\big( \{\e\}^{i-1} \times A \times \{\e\}^{n-i} \big)$;
this preserves joinlessness.

If $f \in n\,{\cal RM}^{\sf fin}_{\rm norm}$:
Then in addition to the properties of $n\,{\cal RM}^{\sf fin}_2$, 
$f({\rm domC}(f)) = {\rm imC}(f)$, i.e., $Q = {\rm imC}(f)$.  
Then in the one-step restriction, $P$ is replaced by $P_{p,i}$, $f$ comes
$f|_R$ where $R = P_{p,i} \, nA^*$, and $Q$ is replaced by 
${\rm imC}(f|_R)$; hence $f|_R(P_{p,i}) = {\rm imC}(f|_R)$. 
So $f|_R$ is normal.
 \ \ \ $\Box$

\begin{lem} \label{LEMhhINVdomC}
 \ For any right-ideal morphism $h$ of $n A^*$:  

\smallskip

\noindent {\small \rm (1)} \ \ \ $h^{-1}(h({\rm domC}(h))) \cdot n A^*$  $=$ 
${\rm Dom}(h)$. 

\smallskip

\noindent {\small \rm (2)} \ \ \ $h^{-1}(h({\rm domC}(h)))$ 
$\,\equiv_{\sf fin}\,$   ${\rm domC}(h)$.
\end{lem}
{\sc Proof.} (1) \ $[\supseteq]$ This follows trivially from the fact that 
$\,{\rm domC}(h) \subseteq$ $h^{-1}(h({\rm domC}(h)))$. 

\noindent $[\subseteq]$ 
For every $x \in h^{-1}(h({\rm domC}(h))) \cdot n A^*\,$ we have
$h(x) \in h(h^{-1}(h({\rm domC}(h)))) = h({\rm domC}(h))\,$ (since 
$h({\rm domC}(h)) \subseteq {\rm Im}(h)$, and for any set 
$S \subseteq {\rm Im}(h)$: $\,h(h^{-1}(S) = S$).
So, $h(x)$ is defined, hence $x \in {\rm Dom}(h)$.
 
\smallskip

\noindent (2) Since ${\rm Dom}(h) = {\rm domC}(h) \cdot n A^*$, item (1)
implies $\, h^{-1}(h({\rm domC}(h))) \cdot n A^*$  $\vartriangle$
${\rm domC}(h) \cdot n A^*$  $=$  $\varnothing$, which is of course finite.
 \ \ \ $\Box$

\begin{pro} \label{restrInRM2Normal} {\bf ({\boldmath 
$\,\equiv_{\sf fin}$-equivalence} to normal).}

\smallskip

\noindent For every $f \in n {\cal RM}^{\sf fin}_0 \,$ 
(i.e.\ ${\rm domC}(f)$ is finite), there exists a restriction $\varphi$ of 
$f$ such that: 

\medskip

\noindent {\small \rm (1)} \ $f \equiv_{\sf fin} \varphi$;    

% and for all $x \in {\rm domC}(f)$ and $z \in {\rm domC}(\varphi)$:
% if $\,x \vee z\,$ exists then $x \le_{\rm init} z$ 

\medskip

\noindent {\small \rm (2)} 
 \ $\varphi \in n\,{\cal RM}^{\sf fin}_{\rm norm} \,$ 
(i.e.\ ${\rm domC}(\varphi)$ and ${\rm imC}(\varphi)$ are finite and 
joinless with $\varphi({\rm domC}(\varphi)) = {\rm imC}(\varphi)$); 

\medskip

\noindent {\small \rm (3)} \ ${\rm maxlen}(\varphi)$ $\,\le\,$ 
     $3 \ {\rm maxlen}(f)$. 
\end{pro}
{\sc Proof.} (1) \& (2).
By applying Prop.\ \ref{initTojoinless} to ${\rm domC}(f)$
we obtain a finite joinless code $C \subseteq {\rm Dom}(f)$ 
such that $\,{\rm domC}(f) \equiv_{\sf fin} C$, and
$\,{\rm maxlen}(C) = {\rm maxlen}({\rm domC}(f)) =_{\rm def} \ell$.
%% Moreover, $C$ is a coordinatewise joinless code, and 
%% $C$ is derivable from $(\arepsilon)^n$ by one-step restrictions. 

Let $\,h = f|_{C \, nA^*}$. Then $\,{\rm domC}(h) = C$ is a finite joinless
code;  hence $h \in n\,{\cal RM}^{\sf fin}_1$. 
Also, since $\,{\rm domC}(f) \equiv_{\sf fin} C$, we have 
$\,h \equiv_{\sf fin} f \ $ (by Lemma \ref{LEMequivMorphRestr}).

Next, by applying Prop.\ \ref{initTojoinless} to $\,h(C) \,$ 
($\, = h({\rm domC}(h))$) we obtain a finite joinless code 
$\,S \subseteq {\rm Im}(h)\,$ such that $\,S \equiv_{\sf fin} h(C)$, 
and $\,{\rm maxlen}(S) = {\rm maxlen}(h(C)) =_{\rm def} L$.
Let $\,D = h^{-1}(S)$. So by Lemma \ref{joinsINimC}(4), $D$ is finite,
and by Lemma \ref{joinsINimC}(3), $D$ is joinless. 

Let $\,\varphi = h|_D$.  Then $\,{\rm domC}(\varphi) = D$, which is finite 
and joinless.  And $\, \varphi({\rm domC}(\varphi))$  $=$  $\varphi(D)$  $=$
$h(D) = h(h^{-1}(S)) = S$. 
Since $S$ is an initial-factor code (more strongly, we saw that $S$ is 
joinless) and $\,S = \varphi(D)$  $=$ $\varphi({\rm domC}(\varphi))$, it 
follows that $\,S = {\rm imC}(\varphi)\,$ (by the definition of {\rm imC}).
Now $\,\varphi({\rm domC}(\varphi)) = {\rm imC}(\varphi)$, hence $\varphi$ 
is normal.

Let us show that $\,D \equiv_{\sf fin} {\rm domC}(f)$. 
By the construction of $S$ in Prop.\ \ref{initTojoinless}(1):  
$S \equiv_{\sf fin} h(C)$, and $S \subseteq h(C) \ nA^*$. 
Moreover, since $S \subseteq h(C) \ nA^*$, we have for all 
$s \in S$ and $y \in h(C)$: if $\,y \vee s\,$ exists, then 
$y \le_{\rm init} s$.  Thus by Lemma \ref{LEMequivInvImage}, 
$\,S \equiv_{\sf fin} h(C)\,$ implies 
$\,h^{-1}(S) \equiv_{\sf fin} h^{-1}(h(C))$. 
By the definition of $D$:  $D = h^{-1}(S)$.
And $h^{-1}(h(C)) = h^{-1}(h({\rm domC}(h)))$  $\equiv_{\sf fin}$ 
${\rm domC}(h)$ \ (by Lemma \ref{LEMhhINVdomC}). Thus, 
$D \equiv_{\sf fin} C$. 
And we saw that $\, C \equiv_{\sf fin} {\rm domC}(f)$.

\medskip

\noindent (3) Let us find an upper bound on 
$\,{\rm maxlen}({\rm domC}(\varphi))$, where $\,{\rm domC}(\varphi) = D$, 
and an upper bound on $\,{\rm maxlen}({\rm imC}(\varphi))$, where 
$\,{\rm imC}(\varphi) = S$.  

\smallskip

\noindent Since $D = h^{-1}(S)$, we have 

\smallskip

 \ \ \  \ \ \ ${\rm maxlen}(D)$ $\le$  
    ${\rm maxlen}(S) + {\rm maxlen}({\rm domC}(h))$.

\smallskip

\noindent By the construction of $S$, 
$\,{\rm maxlen}(S) = {\rm maxlen}(h(C))$. And $h(C) = f(C)$, and 
${\rm domC}(h) = C$.  So, 

\smallskip

 \ \ \  \ \ \ ${\rm maxlen}(D)$ $\le$ 
    ${\rm maxlen}(f(C)) + {\rm maxlen}(C)$.

\smallskip

\noindent And $\,{\rm maxlen}(f(C))$  $\le$ 
${\rm maxlen}(f({\rm domC}(f))) + {\rm maxlen}(C)$.  So, 

\smallskip

 \ \ \  \ \ \ ${\rm maxlen}(D)$ $\le$ 
${\rm maxlen}(f({\rm domC}(f))) +  {\rm maxlen}(C) + {\rm maxlen}(C)$.

\smallskip

\noindent And by the construction of $C$, 
$\,{\rm maxlen}(C) = {\rm maxlen}({\rm domC}(f))$.  

\noindent Thus for ${\rm maxlen}(D)\,$  
($\,= {\rm maxlen}({\rm domC}(\varphi))$) we have:

\medskip

\noindent {\small (3.1)}     
 \ \ \ \  \ \ \ ${\rm maxlen}({\rm domC}(\varphi))$ $\,\le\,$
${\rm maxlen}(f({\rm domC}(f))) +  2 \ {\rm maxlen}({\rm domC}(f))$
$\,\le\,$  $ 3 \ {\rm maxlen}(f)$.

\medskip

\noindent We saw that \ ${\rm maxlen}({\rm imC}(\varphi))$ $=$
${\rm maxlen}(S) \,=\, {\rm maxlen}(h(C))$
$\le\,$  ${\rm maxlen}(f({\rm domC}(f))) + {\rm maxlen}(C)$ 
$\,=\,$  ${\rm maxlen}(f({\rm domC}(f))) + {\rm maxlen}({\rm domC}(f))$.
So,

\medskip

\noindent {\small (3.2)}  
  \ \ \  \ \ \ ${\rm maxlen}({\rm imC}(\varphi))$  $\,\le\,$
${\rm maxlen}(f({\rm domC}(f))) + {\rm maxlen}({\rm domC}(f))$
$\,\le\,$  $ 2 \ {\rm maxlen}(f)$.
 \ \ \ $\Box$

\bigskip

\noindent Recall that the definition of $n M_{k,1}$ was based on
$\, n \,{\cal RM}^{\sf fin}_1 \,$ (Def.\ \ref{DEFnM21}).  Prop.\
\ref{LEMhhINVdomC} now implies:

\begin{cor} \label{nM21RM0RM} {\bf (equivalent definitions of
    {\boldmath $n M_{k,1}$}).}
The following monoids are isomorphic:

\medskip

 \ \ \ $n M_{k,1}$  $=$
$n \,{\cal RM}^{\sf fin}_1 \!/\! \equiv_{\rm fin}\,$,
 \ \ $n \,{\cal RM}^{\sf fin}_0 \!/\! \equiv_{\rm fin}\,$,
 \ \ $n \,{\cal RM}^{\sf fin}_2 \!/\! \equiv_{\rm fin}\,$,
 \ \ $n \,{\cal RM}^{\sf fin}_{\rm norm} /\! \equiv_{\rm fin} \ $.
\end{cor}
{\sc Proof.} By Prop.\ \ref{LEMhhINVdomC}, every element $f$ of
$n \,{\cal RM}^{\sf fin}_0$ is $\,\equiv_{\rm fin}$-equivalent to an element
of $n \,{\cal RM}^{\sf fin}_1$, and to an element of
$n \,{\cal RM}^{\sf fin}_2$, and to a normal element.
 \ \ \ $\Box$

%%%%%%%%%%%%%%%%%%%%%%%%%%%%%%%%%%%%%%%%%%%%%%%%%%%%%%%%%%%%%%%
% Section
%%%%%%%%%%%%%%%%%%%%%%%%%%%%%%%%%%%%%%%%%%%%%%%%%%%%%%%%%%%%%%%
\section{Some algebraic properties of {\boldmath $n M_{k,1}$} }

We show some algebraic properties of $n M_{k,1}$ that are similar to 
properties of $M_{k,1}$ proved in \cite{BiThompsMon, BiThompsMonV3}. The 
proofs are similar too, except for the last cases of Theorem 
\ref{THMnMsimple}.

\bigskip

\noindent Recall that a monoid $M$ is called {\em regular} iff for every
$f \in M$ there exists $f' \in M$ such that $f f' f = f$.

\begin{pro} \label{nM21Regular} {\bf (regularity).} 
 \ The monoid $\,n M_{k,1}$ is regular, for all  $n \ge 1$ and $k \ge 2$.
\end{pro}
{\sc Proof.} For an element of $\,n M_{k,1}$ represented by a right-ideal
morphism $f \in n {\cal RM}^{\sf fin}_{\rm norm}$ (according to Prop.\
\ref{restrInRM2Normal}), we can choose an inverse $f'$ of $f$ with the
following table: $\, {\rm domC}(f') = {\rm imC}(f)\,$
($\,= f({\rm domC}(f))$), $\,{\rm imC}(f') = {\rm domC}(f)$, and for every
$y \in {\rm imC}(f)$: $f'(y)$ is picked arbitrarily in $f^{-1}(y)$. Since
${\rm domC}(f)$ is joinless and ${\rm imC}(f) = f({\rm domC}(f))$ (also
joinless), $f^{-1}({\rm imC}(f)) = {\rm domC}(f)$.
Then $f' \in n {\cal RM}^{\sf fin}_{\rm norm}$, and $f f' f = f$.
 \ \ \ $\Box$

\begin{pro} \label{nM21GroupUnits} {\bf (group of units).} 
 \ The {\em group of units} of $\,n M_{k,1}$ is the Brin-Higman-Thompson 
group $\,n G_{k,1}$, for all  $n \ge 1$ and $k \ge 2$.
\end{pro}
{\sc Proof.} The identity element of $n M_{k,1}$ is the identity function 
${\mathbb 1}$ on $nA^{\omega}$. An element $f \in n M_{k,1}$ belongs to the 
group of units iff there exist $h, k \in n M_{k,1}$ such that 
$h f = f k = {\mathbb 1}$.
Obviously, ${\mathbb 1}$ is a total, surjective and injective function on
$nA^{\omega}$. The fact that $h f(.)$ is total implies that $f$ is total,
and the fact that $f k(.)$ is surjective implies that $f$ is surjective.
And for a total function $f$, the fact that $h f(.)$ is injective implies
that $f$ is injective. Hence, $f \in n G_{k,1}$.
 \ \ \ $\Box$

\bigskip

Recall that in a monoid $M$, two elements $x, y \in M$ are called
${\cal J}$-equivalent, denoted by $x \equiv_{\cal J} y$, iff 
 \ $(\exists \alpha, \beta, \gamma, \delta \in M)[\,x =\beta y \alpha$ 
and $y = \delta x \gamma \,]$. \ (See e.g.\ \cite{Grillet}.)

A monoid $M$ with zero $\mathbb 0$ is called ${\cal J}$-0 {\em simple} iff 
the only ${\cal J}$-equivalence classes of $M$ are 
$\,M \minus \{\mathbb 0\}\,$ and $\{\mathbb 0\}$.

\begin{pro} \label{LEMnM0Jsimple} {\bf ({\boldmath ${\cal J}$-0} 
simplicity).}
 \ The monoid $n M_{k,1}$  is {\em ${\cal J}$-0 simple}, for all  $n \ge 1$
and $k \ge 2$.
\end{pro}
{\sc Proof.} Let $\varphi \in n M_{k,1}$ any non-empty element. Then
$\varphi$ is represented by a right-ideal morphism, also called $\varphi$
here, such that for some $x, y \in nA^*$: $\,\varphi(x) = y$.
Consider $\alpha, \beta \in n M_{k,1}$, where $\alpha$ is given by the
table $\{((\e)^n, \, x)\}$, and $\beta$ is given by the table
$\{(y, \, (\e)^n)\}$. Then $\beta \cdot \varphi \cdot \alpha$
has the table $\{((\e)^n, \, (\e)^n)\}$, which represents
the identity element. Hence $\varphi \ge_{\cal J} {\mathbb 1}$.
So all non-zero elements are ${\cal J}$-equivalent to the identity.
 \ \ \ $\Box$

\bigskip

A monoid $M$ is called {\em congruence-simple} iff the only 
congruences on $M$ are the unavoidable congruences, namely the equality
relation, and the one-class equivalence relation.

\begin{thm} \label{THMnMsimple} {\bf (congruence-simplicity).} 
 \ For all  $n \ge 1$, the monoid $n M_{k,1}$ is {\em congruence-simple}
(i.e., the only congruences are the equality relation, and the one-class
congruence).
\end{thm}
{\bf Proof.} Let $\equiv$ be any congruence on $n M_{k,1}$ that is not
the equality relation. We will show that then the whole monoid is
congruent to the empty element ${\mathbb 0}$ (represented by the empty
right-ideal morphism).

\smallskip

\noindent {\sf Case 0:} \ Assume that there exists $\Phi \in n M_{k,1}$
such that $\Phi \equiv {\mathbb 0}$ and $\Phi \neq {\mathbb 0}$.
Then for all $\alpha, \beta \in n M_{k,1}$ we have obviously
$\alpha \,\Phi\, \beta \equiv {\mathbb 0}$. Moreover, since
$\Phi \neq {\mathbb 0}$ we have \ $n M_{k,1}$
$=$  $\{\alpha \,\Phi\, \beta : \alpha, \beta \in n M_{k,1}\}$,
by 0-$\cal J$-simplicity of $n M_{k,1}$ (Prop.\ \ref{LEMnM0Jsimple}).
Hence in case {\sf 0} all elements of $n M_{k,1}$ are congruent to
${\mathbb 0}$.

\medskip

When we are not in case {\sf 0}, all elements $\psi$ $\in$ 
$n M_{k,1} \minus \{{\mathbb 0}\}$  satisfy $\,\psi \not\equiv {\mathbb 0}$.
For the remainder of the proof we suppose that there exist $\varphi, \psi$ 
$\in$ $n M_{k,1} \minus \{{\mathbb 0}\}$ such that $\varphi \equiv \psi$ and
$\varphi \neq \psi$.
We use the characterization of $\equiv_{\rm fin}$ from Lemma
\ref{LEMequivFinEnds}.

\medskip

Notation: By $(x \mapsto y)$ we denote the element of $n M_{k,1}$ given by
the singleton table $\{(x, y)\}$, where $x, y \in nA^*$.

\medskip

\noindent {\sf Case 1:}
 \ ${\rm Dom}(\varphi) \not\equiv_{\rm fin} {\rm Dom}(\psi)$.

\noindent Then by Lemma \ref{LEMequivFinEnds} we have: There exists
$x \in n A^*$ such that

 \ \ \ $x \cdot nA^* \subseteq {\rm Dom}(\varphi)$, but
 \ ${\rm Dom}(\psi) \,\cap\, x \cdot nA^* = \varnothing$;

\noindent or, vice versa, there exists $x \in n A^*$ such that

 \ \ \ $x \cdot nA^* \subseteq {\rm Dom}(\psi)$, but
 \ ${\rm Dom}(\varphi) \,\cap\, x \cdot nA^* = \varnothing$.

\noindent Let us assume the former (the other case is similar). For
$\beta = (x \mapsto x)$ we have
$\varphi \, \beta(.) = (x \mapsto \varphi(x))$.
We also have $\psi \, \beta(.) = {\mathbb 0}$, since
$x \cdot nA^* \cap {\rm Dom}(\psi) = \varnothing$.
So, $\varphi \, \beta \equiv \psi \, \beta = {\mathbb 0}$, but
$\varphi \, \beta \neq {\mathbb 0}$.
Hence case {\sf 0} can be applied to $\Phi = \varphi \, \beta$, which
implies that the entire monoid $M_{k,1}$ is congruent to ${\mathbb 0}$.

\medskip

\noindent {\sf Case 2.1:}
 \ ${\rm Im}(\varphi) \not\equiv_{\rm fin} {\rm Im}(\psi)$
 \ and \ ${\rm Dom}(\varphi) \equiv_{\rm fin} {\rm Dom}(\psi)$.

\noindent Then there exists $y \in nA^*$ such that
$y \cdot nA^* \subseteq {\rm Im}(\varphi)$, but
$\, {\rm Im}(\psi) \cap y \cdot nA^* = \varnothing$;  or, vice versa,
$y \cdot nA^* \subseteq {\rm Im}(\psi)$, but
$\,{\rm Im}(\varphi) \cap y \cdot nA^* = \varnothing$.  Let us assume the
former (the other case is similar). Let $x \in n A^*$ be such that
$y = \varphi(x)$. Then $\,(y \mapsto y) \circ \varphi \circ (x \mapsto x)$
$\,=\, (x \mapsto y)$.

On the other hand,
$(y \mapsto y) \circ \psi \circ (x \mapsto x) \,=\, {\mathbb 0}$.
Indeed, if $x \cdot nA^* \cap {\rm Dom}(\psi) = \varnothing$ then for all
$v \in A^* : $\, $\psi \circ (x \mapsto x)(xv) \,=\, \psi(xv)$
$\,=\, $ $\varnothing$.
And if $x \cdot nA^* \,\cap\, {\rm Dom}(\psi) \neq  \varnothing$ then
for those $v \in A^*$ such that $xv \in {\rm Dom}(\psi)$ we have
$ \ (y \mapsto y) \circ \psi \circ (x \mapsto x)(xv)$ $\,=\,$
$(y \mapsto y)(\psi(xv)) \,=\, \varnothing$, since
${\rm Im}(\psi) \cap y \cdot nA^* = \varnothing$.  Now case {\sf 0} applies 
to ${\mathbb 0} \neq \Phi$ $=$
$(y \mapsto y) \circ \varphi \circ (x \mapsto x)$  $\equiv {\mathbb 0}$;
hence all elements of $n M_{k,1}$ are congruent to ${\mathbb 0}$.

\smallskip

\noindent {\sf Case 2.2:}
 \ ${\rm Im}(\varphi) \equiv_{\rm fin} {\rm Im}(\psi)$
 \ and \ ${\rm Dom}(\varphi) \equiv_{\rm fin} {\rm Dom}(\psi)$.

\noindent Then after restricting $\varphi$ and $\psi$ to
${\rm Dom}(\varphi) \,\cap\, {\rm Dom}(\psi)$ ($\equiv_{\rm fin}$
${\rm Dom}(\varphi) \equiv_{\rm fin} {\rm Dom}(\psi)$), we have:
$\, {\rm domC}(\varphi) = {\rm domC}(\psi)$, and
 there exist $x \in {\rm domC}(\varphi) = {\rm domC}(\psi)$ and
$y \in {\rm Im}(\varphi)$, $z \in {\rm Im}(\psi)$ such
that $\varphi(x) = y \neq z = \psi(x)$. We have two sub-cases.

\smallskip

\noindent {\sf Case 2.2.1:} \ $y$ and $z$ have no join.

\noindent Then $ \ (y \mapsto y) \circ \varphi \circ (x \mapsto x)$
$=$  $(x \mapsto y)$.  On the other hand,
$\,(y \mapsto y) \circ \psi \circ (x \mapsto x)(x w)$  $\,=\,$
$(y \mapsto y)(z w) \,=\, \varnothing$ \ for all $w \in A^*$ \
(since in this case it is assumed that $y$ and $z$ have no join).  So
$ \ (y \mapsto y) \circ \psi \circ (x \mapsto x) \,=\, {\mathbb 0}$.
Hence case {\sf 0} applies to $ \ {\mathbb 0} \neq\Phi$  $\,=\,$
$(y \mapsto y) \circ \varphi \circ (x \mapsto x)$
$\equiv {\mathbb 0}$.

\smallskip

\noindent {\sf Case 2.2.2:} \ $y \vee z$ exists, and
$y \neq z$.

\noindent Then $y \vee z = zv = yu$ for some
$u, v \in nA^*$, with $v \ne (\e)^n$ or $u \ne (\e)^n$;
let us assume the latter (the other case is similar).
Since $y <_{\rm init} y \vee z$, there exists $v'$ such that $\,z \vee yv'\,$
does not exist (by Lemma \ref{LEMyNzNOJ}).
Now, $\, (yv' \mapsto yv') \circ \varphi \circ (x \mapsto x)(x v')$
$\,=\,$ $(yv' \mapsto yv')(yv')$ $=$  $yv'$.
But for all $s \in nA^*$:
$\, (yv' \mapsto yv') \circ \psi \circ (x \mapsto x)(xs)$   $\,=\,$
$(yv' \mapsto yv')(zs)$  $=$  $\varnothing$, since $\,z \vee yv'\,$
does not exist.
Thus, case {\sf 0} applies to $ \ {\mathbb 0} \neq \Phi = $
$(yv' \mapsto yv') \circ \varphi \circ (x \mapsto x)$
$\equiv$ $\, (yv' \mapsto yv') \circ \psi \circ (x \mapsto x)(xs)$
$\,=\,$   ${\mathbb 0}$.
 \ \ \ $\Box$

%%%%%%%%%%%%%%%%%%%%%%%%%%%%%%%%%%%%%%%%%%%%%%%%%%%%%%%%%%%%%%%%%%%%%%
%%% Section
%%%%%%%%%%%%%%%%%%%%%%%%%%%%%%%%%%%%%%%%%%%%%%%%%%%%%%%%%%%%%%%%%%%%%%
\section{Finite generation of {\boldmath $\,n M_{k,1}\,$} }

\noindent The monoid $n M_{k,1}$ contains the subgroup $n G_{k,1}$ and the
submonoid $\,M_{k,1} \x \{\mathbb{1}_1\}^{n-1}$, where $\mathbb{1}_1$ denotes 
the identity function on $A_k^{\,\omega}$ ($= 1 \, A_k^{\,\omega}$).
It was proved in \cite{BiThompsMon} that $M_{k,1}$ is finitely generated, 
and it was proved in \cite{BinGk1} that $n G_{k,1}$ is finitely generated
(it was previously known that $n G_{2,1}$ is finitely generated
\cite{Brin1, HennigMattucci}).

\begin{pro} \label{nMfactGMG} {\bf (factorization of {\boldmath 
$n M_{k,1}$}).} 

\medskip

\noindent For all $n \ge 2$ and $k \ge 2$, the monoid $\,n M_{k,1}\,$ can 
be factored as 

\medskip

 \ \ \  \ \ \ $n M_{k,1}$ \ $=$  \ $n G_{k,1}$  $\,\cdot\,$ 
$(M_{k,1} \times \{\mathbb{1}_1\}^{n-1})$  $\,\cdot\,$   $n G_{k,1}\,$.

\smallskip

\noindent In other words, for every $f \in n M_{k,1}$ there exist 
$g_1, g_2 \in n G_{,1}$ and $h \in M_{k,1} \x \{\mathbb{1}_1\}^{n-1}$
such that $f = g_2 \, h \, g_1$.
\end{pro}
{\sc Proof.} In order to find such a factorization we use the fact that 
every element of $n M_{k,1}$ can be represented by a function 
$f \in$ $n\,{\cal RM}^{\sf fin}_{\rm norm}$; i.e., $f$ is a right-ideal 
morphism of $n A_k^{\,*}$ such that ${\rm domC}(f)$ and 
${\rm imC}(f)$ are finite joinless codes of $nA_k^{\,*}$ with 
$f({\rm domC}(f)) = {\rm imC}(f)\,$ (Prop.\ \ref{restrInRM2Normal}).

We also use the fact that every maximal finite joinless code in 
$n A_k^{\,*}$ has cardinality $1 + (k-1) \, N$ for some $N \ge 1$; and
that for every $N \ge 1$ there exist maximal prefix codes in $A_k^{\,*}$ of 
cardinality $1 + (k-1) \, N\,$ (see e.g.\ \cite{BerstPerReut},
\cite{BerstelPerrin}, \cite[Coroll.\ 2.14(2)]{BinG}, 
\cite[Lem.\ 9.9(0)]{BiCoNP}). Hence for every maximal finite joinless code 
$C \subseteq n A_k^{\,*}$ there exist a maximal prefix code 
$P \subseteq A_k^{\,*}$ such that $|C| = |P|$.

\smallskip

\noindent $\bullet$ {\sf Case 1.} \ $f$ is {\em total} and {\em surjective},
i.e., ${\rm domC}(f)$ and ${\rm imC}(f)$ ($\, = f({\rm domC}(f))$) are 
maximal finite joinless codes: 

Let $P, Q \subseteq A_k^{\,*}$ be any maximal prefix codes of cardinalities
$\,|P| = |{\rm domC}(f)|$, \ $|Q| = |{\rm imC}(f)|$. 
Let $g_1$ be given by a table which is any bijection from 
$\, {\rm domC}(f)\,$ onto $\,P \x \{\e\}^{n-1}$; and 
let $g_2$ be given by a table which is any bijection from 
$\,Q \x \{\e\}^{n-1}$ onto $\,{\rm imC}(f)$.
Then $g_1, g_2 \in n G_{k,1}$. 
Finally, let $h = g_2^{\, -1} f g_1^{\, -1}(.)$;  \ then 
$h \in M_{k,1} \x \{ \mathbb{1}_1\}^{n-1}$, and $f = g_2 \, h \, g_1$.

\smallskip

\noindent $\bullet$ {\sf Case 2.} \ $f$ is {\em total} but 
{\em not surjective}, i.e., ${\rm domC}(f)$ is a maximal finite joinless 
code, and ${\rm imC}(f)$ ($\, = f({\rm domC}(f))$) is a non-maximal finite 
joinless code:

Since $C = {\rm imC}(f)$ is a non-maximal finite joinless code, it has a 
finite non-empty complementary joinless code (by Coroll.\ 
\ref{CORcomplCode}); let us call it $C'$.

Let $P, Q \subseteq A_k^{\,*}$ be any maximal prefix codes of cardinalities
$\,|P| = |{\rm domC}(f)|$, \ $|Q| = |C \cup C'|$.

Let $g_1 \in n G_{k,1}$ be described by a bijection from
$\,{\rm domC}(f)\,$ onto $\,P \x \{\e\}^{n-1}$, as in Case 1.  

Let $Q_0 \subseteq Q$ be any subset of $Q$ such that $|Q_0| = |C|$, and let 
$Q_0' = Q \minus Q_0$; then $Q_0$ and $Q_0'$ are a pair of complementary
prefix codes in $A_k^{\,*}$, with $|Q_0| = |C|$ and $|Q_0'| = |C'|$.
Let $g_2 \in n G_{k,1}$ be given by a table which is any bijection from $Q$ 
onto $C \cup C'$. 
Then $g_2|_{Q_0 \x \{ \mathbb{1}_1\}^{n-1}}$ is a bijection from 
$Q_0 \x \{ \mathbb{1}_1\}^{n-1}$ onto $C$.
% and $g_2|_{Q_0' \x \{ \mathbb{1}_1\}^{n-1}}$ is a bijection from 
% $Q_0' \x \{ \mathbb{1}_1\}^{n-1}$ onto $C'$.

Then $\,g_2^{\, -1} f g_1^{\, -1}(.)$ $\,=\,$ 
$(g_2|_{Q_0 \x \{ \mathbb{1}_1\}^{n-1}})^{-1} \,f\, g_1^{\, -1}(.)$; and 
$\,h \,=\, g_2^{\, -1} f g_1^{\, -1}(.)$  $\,\in\,$ 
$M_{k,1} \x \{\mathbb{1}_1\}^{n-1}$.

\smallskip

\noindent $\bullet$ {\sf Case 3.} \ $f$ is {\em non-total} but 
{\em surjective},
i.e., ${\rm domC}(f)$ is a non-maximal finite joinless code, and ${\rm imC}(f)$
($\, = f({\rm domC}(f))$) is a maximal finite joinless code:

This is very similar to Case 2.
Since $D = {\rm domC}(f)$ is a non-maximal finite joinless code, it has a
non-empty complementary finite joinless code (by Coroll.\
\ref{CORcomplCode}); let us call it $D'$.

Let $P, Q \subseteq A_k^{\,*}$ be any maximal prefix codes of cardinalities
$\,|P| = |D \cup D'|$, \ $|Q| = |{\rm imC}(f)|$. 

Let $g_2 \in n G_{k,1}$ be as in Case 1.  

Let $P_0 \subseteq P$ be any subset of $P$ such that $|P_0| = |D|$, and let 
$P_0' = P \minus P_0$; then $P_0$ and $P_0'$ are a pair of complementary
prefix codes in $A_k^{\,*}$, with $|P_0| = |D|$ and $|P_0'| = |D'|$. 
Let $g_1 \in n G_{k,1}$ be given by a table which is any bijection from 
$D \cup D'$ onto $P \x \{ \mathbb{1}_1\}^{n-1}$. Then $g_1|_D$ is a 
bijection from $D$ onto $P_0 \x \{ \mathbb{1}_1\}^{n-1}$.
%% , and $g_1|_{D'}$ is a bijection from $D'$ onto 
%% $P_0' \x \{ \mathbb{1}_1\}^{n-1}$.

Then $\,g_2^{\, -1} f g_1^{\, -1}(.)$  $\,=\,$ 
$g_2^{\, -1} \,f\, (g_1|_D)^{-1}(.)$; and 
$\,h = g_2^{\, -1} f g_1^{\, -1}(.)$  $\,\in\,$
$M_{k,1} \x \{\mathbb{1}_1\}^{n-1}$.

\smallskip

\noindent $\bullet$ {\sf Case 4.} \ $f$ is {\em non-total} and
{\em non-surjective}:

Here we combine Cases 2 and 3.
Since $C = {\rm imC}(f)$ is a non-maximal finite joinless code, it has a
non-empty complementary finite joinless code (by Coroll.\
\ref{CORcomplCode}); let us call it $C'$.
Similarly, $D = {\rm domC}(f)$ has a non-empty complementary finite joinless 
code $D'$.

Let $P, Q \subseteq A_k^{\,*}$ be any maximal prefix codes of cardinalities
$\,|P| = |D \cup D'|$, \ $|Q| = |C \cup C'|$.

Let $P_0 \subseteq P$ be any subset of $P$ such that $|P_0| = |D|$, and let 
$P_0' = P \minus P_0$; then $P_0$ and $P_0'$ are a pair of complementary
prefix codes, with $|P_0| = |D|$ and $|P_0'| = |D'|$.
Let $g_1 \in n G_{k,1}$ be given by a table which is any bijection from 
$D \cup D'$ onto $P \x \{ \mathbb{1}_1\}^{n-1}$. Then $g_1|_D$ is a bijection 
from $D$ onto $P_0 \x \{ \mathbb{1}_1\}^{n-1}$.
%, and $g_1|_{D'}$ is a bijection from $D'$ onto 
% $P_0' \x \{ \mathbb{1}_1\}^{n-1}$.

Let $Q_0 \subseteq Q$ be any subset of $Q$ such that $|Q_0| = |C|$, and let 
$Q_0' = Q \minus Q_0$; then $Q_0$ and $Q_0'$ are a pair of complementary
prefix codes, with $|Q_0| = |C|$ and $|Q_0'| = |C'|$.
Let $g_2 \in n G_{k,1}$ be given by a table which is any bijection from 
$Q \x \{ \mathbb{1}_1\}^{n-1}$ onto $C \cup C'$.  Then 
$g_2|_{Q_0 \x \{ \mathbb{1}_1\}^{n-1}}$ is a bijection from 
$Q_0 \x \{ \mathbb{1}_1\}^{n-1}$ onto $C$.
%% $g_2|_{Q_0' \x \{ \mathbb{1}_1\}^{n-1}}$ is a bijection from 
% $Q_0' \x \{ \mathbb{1}_1\}^{n-1}$ onto $C'$.

Then $\,g_2^{\, -1} f g_1^{\, -1}(.)$  $\,=\,$ 
$(g_2|_{Q_0 \x \{ \mathbb{1}_1\}^{n-1}})^{-1} \,f\, (g_1|_D)^{-1}(.)$;
and $\,h \,=\, g_2^{\, -1} f g_1^{\, -1}(.)$    $\,\in\,$
$M_{k,1} \x\{\mathbb{1}_1\}^{n-1}$.
 \ \ \  \ \ \ $\Box$

\bigskip

\noindent {\bf Remark.} The factorization $f = g_2 \, h \, g_1$ in Prop.\ 
\ref{nMfactGMG} is not unique, since the definition of the elements 
$g_2, g_1 \in n G_{k,1}$ allows many arbitrary choices.

\begin{thm} \label{nMfingen}
 \ \ For all $n \ge 1$, $ \ n M_{k,1}\,$ is finitely generated.
\end{thm}
{\sc Proof.} By Prop.\ \ref{nMfactGMG}, every $f \in n M_{k,1}$ has a
factorization $f = g_2 \, h \, g_1$, where $g_1, g_2 \in n G_{k,1}$, 
and $h \in M_{k,1} \x \{ \mathbb{1}_1\}^{n-1}$. Obviously, 
$M_{k,1} \x \{\mathbb{1}_1\}^{n-1}$ is isomorphic to $M_{k,1}$, which is 
finitely generated \cite{BiThompsMon}.
Since $n G_{k,1}$ is also finitely generated \cite{BinGk1}, 
the Theorem follows.
 \ \ \ $\Box$

\bigskip

\noindent {\bf Open question:} Is $n M_{k,1}$ finitely presented?
(This remains open for $M_{k,1}$ too \cite{BiThompsMon}.)

%%%%%%%%%%%%%%%%%%%%%%%%%%%%%%%%%%%%%%%%%%%%%%%%%%%%%%%%%%%%%%%
% Section
%%%%%%%%%%%%%%%%%%%%%%%%%%%%%%%%%%%%%%%%%%%%%%%%%%%%%%%%%%%%%%%
\section{The word problem of {\boldmath $n M_{k,1}$} }

The word problem for a monoid $M$ with a finite generating set $\Gamma$ is
specified as follows:

Input: Two words $u, v \in \Gamma^*$.

Question: $\pi(u) = \pi(v)$, as elements of $M$?  

\noindent Here, for any $w \in \Gamma^*$, $\pi(w) \in M$ denotes the 
product, in $M$, of the generators, as they appear in the word $w$. 
Instead of $\pi(u) = \pi(v)$ we also write $u = v$ in $M$, or 
$u =_M v$.

\begin{lem} \label{LEMcomplexityRel}
 \ Let $M_1$ and $M_2$ be two finitely generated monoids such that 
$M_2 \subseteq M_1$. 

\smallskip

\noindent {\small \rm (1)} \ If the word problem of $M_1$ is in {\sf coNP}
then the word problem of $M_2$ is in {\sf coNP} too.

\smallskip

\noindent {\small \rm (2)} \ If the word problem of $M_2$ is 
{\sf coNP}-hard (with respect to polynomial-time many-one reductions), then 
the word problem of $M_1$ is also {\sf coNP}-hard.
\end{lem}
{\sf Proof.} This is well known; see \cite[Lemma 6.1]{BinGk1} for details.
 \ \ \ $\Box$

\begin{lem} \label{LEMnGk1coNPhard}
 \ The word problem of $n M_{k,1}$ over a finite
generating set is {\sf coNP}-hard, for all $n \ge 2$ and $k \ge 2$.
\end{lem}
{\sc Proof.} In \cite{BinGk1, BinG} it was proved (for all $n \ge 2$, 
$k \ge 2$) that the word problem of the Brin-Higman-Thompson group 
$n G_{k,1}$ over a finite generating set is {\sf coNP}-complete with 
respect to polynomial-time many-one reductions.  Hence the word problem of 
$n M_{k,1}$ over a finite generating set is also {\sf coNP}-hard, by Lemma 
\ref{LEMcomplexityRel}(2).
 \ \ \ $\Box$ 

\bigskip

In the remainder of this section we prove that the word problem of 
$n M_{k,1}$ belongs to {\sf coNP}.

\smallskip

Composition of elements of $n M_{k,1}$, given by tables, works in the 
same way as in $n G_{k,1}$, as described in \cite[Lemma 2.29]{BinG}; in 
outline the proof that the word problem is in {\sf coNP} is the same for 
$n M_{k,1}$ as it was for $n G_{k,1}$, but the intermediary results need 
to be more general.

By Prop.\ \ref{restrInRM2Normal}, every element of $n M_{k,1}$ is represented 
by elements of $n\,{\cal RM}^{\sf fin}_{\rm norm}$; so every element of 
$n M_{k,1}$ has a table of the form $F: P \twoheadrightarrow Q$ where $P$ 
and $Q$ are finite joinless codes, and $F(P) = Q$.

\medskip

\noindent {\bf Notation.} We abbreviate ${\rm maxlen}(.)$ by $\ell(.)$.

\begin{lem} \label{LEMnMcomposition} {\bf (composition in
{\boldmath $n M_{k,1}$} based on tables).}
 \ Let $f_j$: $P_j \twoheadrightarrow Q_j$ be tables of functions in 
$\,n{\cal RM}^{\sf fin}_{\rm norm}$ representing elements of $n M_{k,1}$;
so $P_j, Q_j$ are finite joinless codes, and $f_j(P_j) = Q_j\,$ (for 
$j = 1, 2$).
Then the {\em composite} $f_2 \circ f_1(.)$ is represented by the table

\medskip

 \ \ \  \ \ \ $(f_2 \circ f_1)|_P: \, P \twoheadrightarrow Q$, \ \ where

\medskip

 \ \ \  \ \ \ $P = f_1^{-1}(P_2 \vee Q_1)$,

\medskip

 \ \ \  \ \ \ $Q = f_2(P_2 \vee Q_1)$.
\end{lem}
{\sc Proof.} In general, for any partial functions $f_2, f_1$, we have:
$\, {\rm Dom}(f_2 \circ f_1) = f_1^{-1}({\rm Dom}(f_2) \cap {\rm Im}(f_1))$
$=$ $f_1^{-1}({\rm Dom}(f_2))$, 
and $\, {\rm Im}(f_2 \circ f_1) = f_2({\rm Dom}(f_2) \cap {\rm Im}(f_1))$
$=$   $f_2({\rm Im}(f_1))$.
%Obviously, $f_2 \circ f_1 = (f_2 \circ f_1)|_{{\rm Dom}(f_2 \circ f_1)}$.

When $f_2, f_1$ are given by tables as above,
${\rm Dom}(f_2) \cap {\rm Im}(f_1) = (P_2 \vee Q_1) \cdot nA^{\omega}\,$ 
(by Lemma \ref{SetJoin}). And $f_1^{-1}(P_2 \vee Q_1)\,$ is a finite joinless 
code (by Lemma \ref{joinsINimC}). Moreover, 
$f_1^{-1}((P_2 \vee Q_1) \cdot nA^{\omega})$  $=$ 
$f_1^{-1}(P_2 \vee Q_1) \cdot nA^{\omega}$, by Lemma \ref{LEMinvimage}. The
latter Lemma has the hypothesis, called condition $(\star)$, that for all 
$q \in P_2 \vee Q_1$ and all $y \in f_1({\rm domC}(f_1))\,$ ($\,= Q_1$): 
if $q \vee y$ exists then $y \le_{\rm init} q$. This hypothesis holds by the
connection between $P_2 \vee Q_1$ and $Q_1$.  This shows that 
$\,{\rm Dom}(f_2 \circ f_1) = f_1^{-1}(P_2 \vee Q_1) \cdot nA^{\omega}$,
and $\,{\rm domC}(f_2 \circ f_1) = f_1^{-1}(P_2 \vee Q_1)$.

Finally, $f_2 f_1(f_1^{-1}(P_2 \vee Q_1) \cdot nA^{\omega})$  $=$ 
$f_2 f_1 f_1^{-1}((P_2 \vee Q_1) \cdot nA^{\omega})$  $=$ 
$f_2((P_2 \vee Q_1) \cdot nA^{\omega})$  $=$
$f_2(P_2 \vee Q_1) \cdot nA^{\omega}$; the first equality holds by Lemma
\ref{LEMinvimage}, as we saw; the last equality holds because
$f_2$ is a right-ideal morphism and $P_2 \vee Q_1 \subseteq {\rm Dom}(f_2)$.
This shows that 
$\,{\rm Im}(f_2 \circ f_1) = f_2(P_2 \vee Q_1) \cdot nA^{\omega}$, and 
$\,f_2 \circ f_1({\rm domC}(f_2 \circ f_1)) = f_2(P_2 \vee Q_1)$.

We also have $\, f_2 f_1(P) = f_2f_1f_1^{-1}(P_2 \vee Q_1)$ 
$=$  $f_2(P_2 \vee Q_1)$  $=$ $Q$.
Hence, $f_2 \circ f_1$ is given by the table described in this Lemma.
 \ \ \ $\Box$

\medskip

\noindent {\bf Remark.} In the table for $f_2 \circ f_1$ in Lemma 
\ref{LEMnMcomposition}, $\,Q = f_2 f_1(P)\,$ is not always joinless, since
$\,n{\cal RM}^{\sf fin}_{\rm norm}$ is not closed under composition.

\begin{pro} \label{LengthFormula} {\bf (length formula).} 
 \ For all $f_2, f_1 \in n {\cal RM}^{\sf fin}_{\rm norm}$:
 \ \ $\ell(f_2 \circ f_1) \, \le \, \ell(f_2) + \ell(f_1)$.
\end{pro}
{\sc Proof.} In outline the proof is similar to the one for
\cite[Prop.\ 3.2]{BinG}. 

Let $f_j|_{P_j}: P_j \twoheadrightarrow Q_j$ be a table for $f_j$ ($j=1,2$),
where $P_j$ and $Q_j$ are finite joinless codes, and $f_j(P_j) = Q_j$.
Recall the table for $f_2 \circ f_1$, given in Lemma \ref{LEMnMcomposition}.

\noindent We have:

\smallskip

\noindent (L1) \hspace{1.3in}
$\ell(P_2 \vee Q_1) = \max\{ \ell(P_2), \ \ell(Q_1)\}$,

\medskip

\noindent since for every $p = (p_1, \, \ldots, p_n) \in P_2$,
 \ $q = (q_1, \, \ldots, q_n) \in Q_1$, and $i \in \{1, \, \ldots, n\}$:
 \ $|(p \vee q)_i| = \max\{|p_i|, \, |q_i|\} \,$ (by Lemma
\ref{LEMjoinCharacteriz}).

\noindent We have:

\smallskip

\noindent (L2) \hspace{1.3in} $\ell(f_2(P_2 \vee Q_1))$   $\,\le \,$
$\ell(Q_2) + \ell(Q_1)$  $\,\le\,$   $\ell(f_2) + \ell(f_1)$.

\medskip

\noindent Indeed, $(p \vee q)_i = \max_{\le_{\rm pref}}\{p_i, q_i\}$,
for every $p \in P_2$, $q \in Q_1$, and $i \in \{1, \, \ldots, n\}$.
By Lemma \ref{SetJoin}, $p \vee q \in P_2 \vee Q_1 \subseteq P_2 \cdot nA^*$
$=$  ${\rm Dom}(f_2)$.
Since $p$ is an initial factor of $p \vee q$ there exists $u \in nA^*$ such
that $p u = p \vee q$.
Since $(p \vee q)_i = \max_{\le_{\rm pref}}\{p_i, q_i\}$, the following
holds: $\, u_i = \e$ when $(p \vee q)_i = p_i$; and
$u_i$ is a suffix of $q_i$ when $(p \vee q)_i = q_i$.
Hence, $\ell(u) \le \ell(q)$.
Now, $f_2(p \vee q) = f_2(p) \, u$, where $f_2(p) \in Q_2$ (since
$p \in P_2$). And $q \in Q_1$. Hence 
$\ell(f_2(p \vee q)) \le \ell(f_2(p)) + \ell(u)$
$\le$ $\ell(Q_2) + \ell(Q_1)$.

\noindent We also have:

\smallskip

\noindent (L3) \hspace{1.3in} $\ell(f_1^{-1}(P_2 \vee Q_1))$  $\,\le\,$
$\ell(P_1) + \ell(P_2)$ $\, \le \,$ $\ell(f_2) + \ell(f_1)$.

\medskip

\noindent Indeed, $z \in f_1^{-1}(P_2 \vee Q_1)$ iff $z \in P_1 \, nA^*$
$=$  ${\rm Dom}(f_1)$ and $f_1(z) \in P_2 \vee Q_1$. 
Since $z \in P_1 \, nA^*$ we have $z = p^{(1)} v$ for some $p^{(1)}\in P_1$ 
and $v \in nA^*$, hence $f_1(z) = f_1(p^{(1)}) \, v$; and 
$f_1(p^{(1)}) = q^{(1)}$ for some $q^{(1)}\in Q_1$. 
Since $f_1(z) = q^{(1)} v \in P_2 \vee Q_1$, we have $q^{(1)} v = p \vee q$ 
for some $p \in P_2$ and $q \in Q_1$; since $Q_1$ is joinless, we have 
$q = q^{(1)}$.
Since $(p \vee q^{(1)})_i = \max_{\le_{\rm pref}}\{p_i, q^{(1)}_i\}$, the 
following holds: $\, v_i = \e$ if $(p \vee q^{(1)})_i = q^{(1)}_i$;
and $v_i$ is a suffix of $p_i$ if $(p \vee q^{(1)})_i = p_i$.  
In any case, $\ell(v) \le \ell(p) \le \ell(P_2)$.

So for every $z \in f_1^{-1}(P_2 \vee Q_1)$: 
$\, \ell(z) = \ell(p^{(1)} v)$  $\le$  $\ell(p^{(1)}) + \ell(v)$  $\le$
$\ell(P_1) + \ell(P_2)$.

\smallskip

\noindent Finally, since $\, \ell(f_2 \circ f_1)$ $\, = \,$
$\max\{\ell(f_1^{-1}(P_2 \vee Q_1))$, $\ell(f_2(P_2 \vee Q_1))\}$, we 
conclude:

\medskip

\noindent (L4) \hspace{.5in}
$\, \ell(f_2 \circ f_1) \, \le \, \max\{\ell(Q_2) + \ell(Q_1), \ $
$\ell(P_1) + \ell(P_2)\}$ $\, \le \,$ $\ell(f_2) + \ell(f_1)$.

\smallskip

\noindent $\Box$

\smallskip

\begin{pro} \label{LengthPROD}\!\!. 

\smallskip

\noindent {\small \bf (1)} \ Let 
$f_t, \, \ldots, f_1 \in n {\cal RM}^{\sf fin}_{\rm norm}$, and suppose 
$c \in {\mathbb N}$ is such that $\ell(f_j) \le c\,$ (for 
$j = 1, \, \ldots, t$).
Then there exists $\varphi \in n {\cal RM}^{\sf fin}_{\rm norm}$ such that:

\smallskip

\noindent $\bullet$ \hspace{0.2in} 
$\varphi \,\equiv_{\rm fin}\, f_t \circ \, \ldots \, \circ f_1$ ;

\smallskip

\noindent $\bullet$ \hspace{0.2in}
${\rm Dom}(\varphi)$  $\,\subseteq\,$ 
${\rm Dom}(f_t \circ \, \ldots \, \circ f_1)$ ;

\smallskip

\noindent $\bullet$ \hspace{0.2in}
${\rm maxlen}(\varphi) \,\le\, 6 c t^3$. 

\medskip

\noindent {\small \bf (2)} \ Hence for all $\,\ell \ge 6 c t^3$:  \ \   
$\varphi|_{nA^{\ell}} \ = \ (f_t \circ \,\ldots\, \circ f_1)|_{nA^{\ell}}$.
\end{pro}
{\sc Proof.} (1) If $t$ is not a power of 2, we replace $t$ by $T$   $=$
$2^{\lceil \log_2 t \rceil}$, which is the nearest power of 2 larger than 
$t$; then $\,t \le T < 2t$.
By adding $\,T - t\,$ identity functions to the list of functions we obtain a 
product $f_T \circ \,\ldots\, \circ f_{t+1}$  $\circ$ 
$f_t \circ \,\ldots\, \circ f_1$, where 
$f_T = \,\ldots\, = f_{t+1} = {\mathbb 1}$, where $T$ is a power of 2.

We will prove the existence of $\varphi$ by induction on $\log_2 T$ (which 
is an integer for powers of 2). 
For the length we will actually prove a tighter relation: 
$\,\ell(\varphi) \,\le\,$  $c \ 6^{\log_2 T}\,$ ($\,= c \, T^{\log_2 6}$). 

For $\log_2 T = 0$ (i.e., $T = 1$), the result is obvious.
Inductively, when $\log_2 T \ge 1$ (i.e., $T \ge 2$), we 
subdivide the list of functions into pairs 

\smallskip

 \ \ \ $(f_T \circ f_{T-1})$  $\circ$  $\ldots\,$   $\circ$  
$(f_{2j} \circ f_{2j-1})$  $\circ$ $\,\ldots\,$  $\circ$  $(f_2 \circ f_1)$.

\smallskip

\noindent Let $F_j = f_{2j} \circ f_{2j-1}$.  
By Prop.\ \ref{LengthFormula}, 

\smallskip

$\ell(F_j) \,\le\, \ell(f_{2j}) + \ell(f_{2j-1}) \le 2c\,$. 

\smallskip

\noindent  Since $F_j$ does not necessarily belong to 
$n {\cal RM}^{\sf fin}_{\rm norm}$, we apply Prop.\ \ref{restrInRM2Normal} 
to obtain an $\,\equiv_{\rm fin}$-equivalent element $h_j \in$ 
$n {\cal RM}^{\sf fin}_{\rm norm}$; this is needed for the inductive step.
By Prop.\ \ref{restrInRM2Normal}, 
$\,F_j \equiv_{\rm fin} h_j$ for some  
$h_j \in n {\cal RM}^{\sf fin}_{\rm norm}$ such that
$\,\ell(h_j) \,\le\, 3 \ \ell(F_j)$  $\,\le\, 6 c$.
This yields a product of $\frac{T}{2}$ functions $h_j \in$ 
$n {\cal RM}^{\sf fin}_{\rm norm}$, with $\ell(h_j) \,\le\, 6c\,$ (for 
$j = 1, \ldots, \frac{T}{2}$). 
Since $\log_2 \frac{T}{2} = \log_2 T - 1$, induction now implies that there 
exists $\varphi \in n {\cal RM}^{\sf fin}_{\rm norm}$ such that:
$ \ \varphi \,\equiv_{\rm fin}\, h_{\frac{T}{2}} \circ \, \ldots \,$ 
$\circ$  $h_1$, 
 \ \ ${\rm Dom}(\varphi)$  $\,\subseteq\,$
${\rm Dom}(h_{\frac{T}{2}} \circ \, \ldots \, \circ h_1)$, and
 \ $\ell(\varphi) \,\le\, 6c \ 6^{\log_2 \frac{T}{2}}\,$ 
($\, = c \ 6^{\log_2 T}$).
Since $\,h_{\frac{T}{2}} \circ \, \ldots \, \circ$  $h_1$   $=$
$f_T \circ \,\ldots\, \circ f_1$, we have the claimed result.

Finally, for the length we have: \ $c \, 6^{\log_2 T}$  $=$
$c \, T^{\log_2 6}$  $<$  $c \, (2t)^{\log_2 6} = 6 \, c \, t^{\log_2 6}$
$<$  $6 \, c \, t^3$; here we used the fact that $T < 2t$, and
$\log_2 6 < 3$.     %% $\log_2 6 \simeq 2.584962500721156$ 
 
\smallskip

\noindent (2) This follows now from Lemma \ref{LEMequivMorphBound}.
 \ \ \ $\Box$

\bigskip

\noindent {\bf Definition of {\sf coNP} and {\sf NP}:} \ We use the 
logic-based version of the definitions (see e.g.\ \cite{HandbookA},
\cite{BinG}).
A set $S \subseteq \Gamma^* \x \Gamma^*$ is in {\sf coNP} \ iff \ there 
exists $m \ge 1$, a two-variable predicate 
$\,R(.,.) \,\subseteq\, mA^* \times (\Gamma^* \x \Gamma^*)$, and a 
polynomial $p(.)$, such that

\smallskip

\noindent {\small (1)} \ \  \ \ \ $R \in {\sf P}$ \ (i.e., membership in 
$R$ is decidable in deterministic polynomial time);

\smallskip

\noindent {\small (2)} \ \  \ \ \ $S$  $\,=\,$
$\{(u,v) \in \Gamma^* \x \Gamma^* \,: \ $
$(\forall x \in m A^{\le p(|u| + |v|)}) \ R(x, (u,v)) \, \}$.

\smallskip

\noindent (The definition of {\sf NP} is similar, replacing $\forall$ by
$\exists$.)

\begin{lem} \label{nM21inCoNP}
 \ The word problem of $n M_{k,1}$ over any finite generating set belongs to
{\sf coNP}.
\end{lem}
{\sc Proof.} Let $\Gamma$ be any finite generating set of $n M_{k,1}$.
Every $\gamma \in \Gamma$ has a finite table $F_{\gamma}$: 
$P_{\gamma} \twoheadrightarrow Q_{\gamma}$.
For any $w \in \Gamma^*$, let $f_w \in n {\cal RM}^{\sf fin}_1$ be the 
function obtained by composing the generators in $w$, and let 
$\varphi_w \in n {\cal RM}^{\sf fin}_{\rm norm}$ be the  
$\,\equiv_{\rm fin}$-equivalent normal element, according to Prop.\
\ref{LengthPROD}: $\,f_w \equiv_{\rm fin} \varphi_w$, and
$\, \ell(f_w) \le \ell(\varphi_w) \le c_{_{\Gamma}} \, |w|^3$, where
$\, c_{_{\Gamma}} = 6 \ \max\{\ell(\gamma) : \gamma \in \Gamma\}$.  So
$c_{_{\Gamma}}$ is a known constant, determined by $\Gamma$. And 
$f_w|_{nA^m} = \varphi_w|_{nA^m}$ for all $m \ge c_{_{\Gamma}} \, |w|^3$, 
by Lemma \ref{LEMequivMorphBound}.
Then for the word problem we have:

\smallskip

\hspace{0.4in}  $u = v\,$ in $nM_{k,1}$ \ \ iff 
 \ \ $\, f_u|_{nA^{L(u,v)}} = f_v|_{nA^{L(u,v)}}\,$ in 
$n {\cal RM}^{\sf fin}_{\rm norm}$,

\smallskip 

\noindent where 
$\,L(u,v) = \max\{c_{_{\Gamma}} \, |u|^3, \ c_{_{\Gamma}} \, |v|^3\}$.
Thus, we obtain the following {\sf coNP}-formula for the word problem:

\smallskip

\hspace{0.4in}  $u = v\,$ in $n M_{k,1}$ \ \ iff
 \ \ $(\forall x \in n A^{L(u,v)})$ $[ \, f_u(x) = f_v(x) \,]$.

\medskip

\noindent We still need to show that the predicate $R(.,.)$, defined by

\smallskip

$R(x, (u,v))$ \ $\Leftrightarrow$
 \ $[ \, (\forall i \in \{1,\ldots,n\})[|x_i| = L(u,v)]$  $\,\Rightarrow\,$
$f_u(x) = f_v(x) \,]$,

\smallskip

\noindent belongs to {\sf P}. \ I.e., we want a deterministic 
polynomial-time algorithm that on input $x \in n A^{L(u,v)} \, $ and
$(u,v) \in \Gamma^* \x \Gamma^*$, checks whether $f_u(x) = f_v(x)$.
To do this we apply, to $x \in n A^{L(u,v)}$, the tables of the
generators $\gamma_j, \delta_k \in \Gamma$ that appear in $u =$ 
$\gamma_t \,\ldots\, \gamma_1$, and in $v = \delta_s \,\ldots\, \delta_1$.  
For $u$ we compute 

\smallskip

$\, x$ \ $\longmapsto$ \ $\gamma_1(x) = y^{(1)}$  
 \ $\longmapsto$ \ $\gamma_2(y^{(1)}) = y^{(2)}$ \ $\longmapsto$
 \ \ $\ldots$ \ \ $\longmapsto$ \ $\gamma_t(y^{(t-1)}) = y^{(t)} = f_u(x)$.

\smallskip

\noindent Since $x \in n A^{L(u,v)}$ $\subseteq$ ${\rm Dom}(f_u)$, every
$y^{(j)}$ is defined.

\smallskip

In the complexity analysis we use the following notation: For $z \in nA^*$,
let $\,|z| = \sum_{i=1}^n |z_i| \ $ (i.e., the sum of the lengths of the 
coordinate strings).
 
By repeatedly using Prop.\ \ref{LengthFormula} we conclude that
$\, |y^{(j)}| \,\le\, |x| + c_{_{\Gamma}}\, j$  $\,\le\,$ 
$n \,L(u,v) + c_{_{\Gamma}}\, |u|$  $\,\le\,$   $2n\, L(u,v)$. 
We proceed in a similar way with $v$ on input $x$. 
After these two computations we check whether $f_u(x) = f_v(x)$.
The application of the table of $\gamma_j$ to $y_{j-1}$ takes time
proportional to $|y_{j-1}| \,$ (for $j = 1, \, \ldots, t$).
So, the time complexity of verifying whether $x$ and $(u,v)$ satisfy the
predicate $R$ is, up to a constant factor, bounded by 
$\,\le\,$  $|x| + \sum_{j=1}^t |y_j|$ $\,\le\,$
$n \,L(u,v) + (|u| + |v|) \cdot n \,L(u,v)$.  And we saw that 
$L(u,v) \,\le\, $ $c_{_{\Gamma}} \, (|u|^3 + |v|^3)$. 
Hence the time-complexity of the predicate $R$ is polynomially bounded
in terms of $|u| + |v|$.
 \ \ \  \ \ \ $\Box$

\begin{thm} \label{THMwpCoNP}
 \ For all $n \ge 2$ and $k \ge 2$, the word problem of $n M_{k,1}$ over 
a finite generating set is {\sf coNP}-complete (with respect to 
polynomial-time many-one reductions).
\end{thm}
{\sc Proof.} By Lemma \ref{nM21inCoNP}, the problem is in {\sf coNP}.
Moreover, by Lemma \ref{LEMnGk1coNPhard}, it is {\sf coNP}-hard.
 \ \ \ $\Box$

\bigskip

\bigskip

%%%%%%%%%%%%%%%%%%%%%%%%%%%%%%%%%%%%%%%%%%%%%%%%%%%%%%%%%%%%%%%
% Section
%%%%%%%%%%%%%%%%%%%%%%%%%%%%%%%%%%%%%%%%%%%%%%%%%%%%%%%%%%%%%%%
\section{A connection between {\boldmath $2 M_{2,1}$} and acyclic circuits}
%% \noindent {\Large \bf Remark: 

In \cite{BiThompsMonV3}, \cite{BiThompsMon} it was proved that $M_{2,1}$ 
is finitely generated; let $\Gamma$ be a finite generating set of $M_{2,1}$.
In \cite[Section 2]{Bi1wPerm}, \cite{BiThompsMonV3}, \cite{BiThompsMon} it 
was proved that for every acyclic circuit $C$ of size $|C|$ there exists a 
word $w_C \in (\Gamma \cup \tau)^*$ such that: 

\smallskip 

\noindent $\bullet$ \ \  For every circuit $C$ and every 
$x \in \{0,1\}^{\ge N_C}$ (where $N_C$ is the input-size of $C$): 
 \ $C(x) = w_C(x)\,$.

\smallskip

\noindent $\bullet$ \ \ There exists a polynomial $p(.)$ such that
for  every circuit $C$: \ \ $|C| \,\le\, |w_C| \,\le\, p(|C|)\,$.

\smallskip

\noindent Here, $\tau = \{ \tau_{j,j+1}:  j \ge 1\}$ is the set of 
bit-position transpositions; 
$\,\tau_{j,j+1}(x_1 \ldots x_{j-1} x_j x_{j+1} \,v)$ $\,=\,$ 
$x_1 \ldots x_{j-1} x_{j+1} x_j \,v$, \ for any $x_1,$ $\ldots,$ $x_{j-1},$ 
$x_j,$  $x_{j+1} \in \{0,1\}$ and $v \in \{0,1\}^*$. 
By $|w_C|$ we denote the size of the word $w_C$ over $\Gamma \cup \tau$, 
where $|\gamma| = 1$ for every $\gamma \in \Gamma$, and 
$\,|\tau_{j,j+1}| = j+1$.  
(See \cite{Bi1wPerm,BiThompsMonV3,BiThompsMon} for details.)

Hence, words over the infinite generating set $\Gamma \cup \tau$ of 
$M_{2,1}$ can represent any acyclic circuits, with at most polynomial 
distortion. 

\medskip

\noindent Since $M_{2,1}$ is finitely generated, acyclic circuits can also 
be represented by words over a finite generating set. But over a finite 
generating set, this representation has exponential distortion; in 
fact, the minimum-length words over $\Gamma$ that represent the 
transpositions $\tau_{j,j+1}$ grow exponentially in length with $j$. 

\medskip

\noindent In \cite{BinG} it was proved that $\tau_{j,j+1} \x {\mathbb 1}_1$ 
can be expressed by a word of length $\le 2j$ over 
$\tau_{1,2} \x {\mathbb 1}_1$ and $\sigma$; here, ${\mathbb 1}_1$ denotes 
the identity function on $\{0,1\}^*$, and $\sigma \in 2 M_{2,1}$ is the 
shift.  In combination with the above this yields:

\medskip

\noindent  {\it 
The finitely generated monoid $2 M_{2,1}$ provides a way to represent
acyclic digital circuits by words over the finite generating set of
$2 M_{2,1}$.  This description preserves sizes up to a polynomial. 
The word problem of $2 M_{2,1}$ over a finite generating set is 
polynomial-time equivalent to the equivalence problem for acyclic circuits 
(and both problems are {\sf coNP}-complete).
}

%%%%%%%%%%%%%%%%%%%%%%%%%%%%%%%%%%%%%%%%%%%%%%%%%%%%%%%%%%%%%%%
%%%%%%%%%%%%%%%%%%%%%%%%%%%%%%%%%%%%%%%%%%%%%%%%%%%%%%%%%%%%%%%

\bigskip

\bigskip

%%%%%%%%%%%%%%%%%%%%%%%%%%%%%%%%%%%%%%%%%%%%%%%%%%%%%%%%%%%%%%%%
%%%%%%%%%%%%%%%%%%%%%%%%%%%%%%%%%%%%%%%%%%%%%%%%%%%%%%%%%%%%%%%%

%% \small

%%%%%%%%%%%%%%%%%%%%%%%%%%%%%%%%%%%%%%%%%%%%%%%%%%%%

\medskip

\noindent {\scriptsize  birget@camden.rutgers.edu }

\end{document}